%% file: cpd_anton_aug26.tex
\documentclass[a4paper,10pt]{article}

\usepackage[T1]{fontenc}
\usepackage{comment,todonotes}
\usepackage{booktabs}
\usepackage{longtable}
\usepackage{hyperref}
\usepackage{amsfonts,amsmath,amsthm}

\usepackage[margin=25mm]{geometry}

\usepackage[backend=bibtex]{biblatex}
\usepackage{tikz}
\usepackage{tikz-3dplot}
\usetikzlibrary{arrows.meta,calc,positioning}

\title{From Continuum-Kinematics-Inspired Peridynamics to Simplicial Hyperelasticity}
\author{Anton Evgrafov\thanks{%
Department of Mathematical Sciences, University of Aalborg, Denmark, anev@math.aau.dk}, 
Pablo Pedregal\thanks{%
Department of Mathematics, University of Castilla-la Mancha, Spain, Pablo.Pedregal@uclm.es}}
\date{\today}

\newcommand{\R}{\mathbb{R}}
\newcommand{\bx}{\boldsymbol{x}}

\newcommand{\bu}{\boldsymbol{u}}
\newcommand{\bv}{\boldsymbol{v}}

\newcommand{\bD}{\boldsymbol{D}}
\renewcommand{\O}{\Omega}
\newcommand{\mres}{\mathbin{\llcorner}}
\newtheorem{theorem}{Theorem}[section]
\newtheorem{corollary}[theorem]{Corollary}
\newtheorem{definition}{Definition}[section]
\newtheorem{lemma}[theorem]{Lemma} 
\newtheorem{proposition}[theorem]{Proposition}

\newtheorem{remark}[theorem]{Remark}
\providecommand{\keywords}[1]{% 
\par\smallskip \noindent\textbf{Keywords:} #1\par } 
\providecommand{\subjclass}[1]{% 
\par\smallskip \noindent\textbf{2020 Mathematics Subject Classification:} #1\par }

\begin{document}

\maketitle

\begin{abstract} 
We develop a variational and functional-analytic framework for continuum-kinematics-inspired peridynamics with simultaneous nonlocal interactions among multiple points.
We encode the finite-difference data of each interacting simplex in a simplicial deformation gradient.
For nondegenerate simplices in three dimensions, the gradient, its cofactor, and its determinant encode the associated continuum-kinematics-inspired peridynamic line, area, and volume changes, respectively. 
This unified representation leads to a class of simplicial hyperelastic energies with general, possibly nonadditive, constitutive dependence on the affine deformation of an interacting simplex, thereby extending the original continuum-kinematics-inspired peridynamic framework.

To describe interactions occurring at different scales and over different simplex geometries, we introduce radial- and shape-weighted nonlocal spaces. 
Their multipoint seminorms control two-point energies induced by the edge marginals of the interaction measure. 
This comparison yields nonlocal Poincar{\'e} inequalities and fractional Sobolev embeddings.

We establish existence of equilibrium configurations by two complementary mechanisms. 
Particle-wise averaged convexity yields weak lower semicontinuity of the associated energy functionals. 
Under stronger assumptions on the edge marginals, fractional compactness gives existence for nonconvex energy densities. 
These results provide a mathematical foundation for simplicial and continuum-kinematics-inspired peridynamic models.
\end{abstract}


\section{Introduction and motivation}
\label{sec:intro}

Although nonlocal models of solid mechanics arguably predate their more conventional differential statements~\cite{Darrigol2002, Navier1821}, recent years have witnessed a renewed interest in nonlocal approaches.
In particular, models involving pairwise nonlocal interactions between particles in the original spirit of Navier, such as peridynamics~\cite{Silling2000} and related theories of nonlocal elasticity~\cite{Kunin1975}, have attracted considerable attention.
These models naturally give rise to variational statements, in which pairwise interactions correspond to double integrals~\cite{Du2019, mengesha2015variational}.
For a broader variational perspective on nonlocality and the associated parametrized-measure framework, we refer to~\cite{MR4330491,Pedregal1997}.

Bond-based models are subject to well-known constitutive restrictions and, in general, do not recover arbitrary classical hyperelastic energies in the vanishing-horizon limit~\cite{MR4150441}.
In the mechanical engineering community, the constitutive restrictions of bond-based peridynamics  motivated the development of state-based peridynamics, in which the response at a material point may depend collectively on the family of incident bonds~\cite{SillingEtAl2007}.
A complementary line of research develops nonlocal analogues of classical differential operators, including gradients, divergences, and integration-by-parts formulas, with applications to peridynamics and nonlocal diffusion; see, for example, \cite{DuEtAlNonlocalCalculus,AlaliLiuGunzburger2015}. 
Related variational models replace the classical distributional gradient by suitably constructed nonlocal integral operator or fractional gradients~\cite{MR4179861, MR4626320,MR4651657,MR4881606,MR3954360}.
This line of attack has been quite successful in the recent mathematical literature and has led to increasingly complete lower-semicontinuity and localization theories for functionals involving nonlocal gradients, see for example~\cite{MR4767360,MR4330183}.
In several such formulations, the relevant structural conditions are similar to those known from the local theory, that is, to quasiconvexity and polyconvexity.

In this work, we study a different class of nonlocalities.
We will consider models allowing for \emph{simultaneous multi-point} nonlocal interactions involving up to \(N+1\) points in \(N\) spatial dimensions.
Such interactions appear, for example, in continuum-kinematics-inspired peridynamics (CPD), see~\cite{javili2019continuum,javili2020computational,steinmann2025foundations}, which serves as our primary inspiration for this study.
Our objective is to propose a suitable mathematical framework, which is general enough to encompass a range of nonlocal models of continuum mechanics, including the original bond-based peridynamics~\cite{Silling2000} and CPD~\cite{javili2019continuum}.
Mathematically, our framework builds upon Pedregal's study of nonlocal variational problems involving multiple integrals~\cite{pedregal_special_2024}.

A central object of our framework is that of a simplicial gradient.
Specifically, let \(\bu: \R^N \to \R^d\) be a mapping.
For mechanical applications the case \(N=d \in \{2,3\}\) will be most pertinent, and the role of \(\bu\) will be to describe displacements of points.
Furthermore, let us fix \(m\in \{1,\ldots,N\}\) and consider a simplex 
\(\operatorname{conv}\{ \bx_0, \bx_1, \ldots, \bx_{m} \}\) \hbox{ with 
vertices }\(\bx_i\in\R^N\), \(i=0,\dots,m\).
Its edges based at \(\bx_0\) are given by vectors \(\boldsymbol{h}_i = \bx_i - \bx_0\), \(i=1,\ldots,m\), which we will assemble into a
matrix \(\boldsymbol{H} = [\boldsymbol{h}_1,\ldots,\boldsymbol{h}_m]\in \R^{N\times m}\).
If our goal was to reconstruct some information about the first derivatives of \(\bu\) based on its values at the simplex vertices, we could
look at first-order approximations of \(\bu\) at \(\bx_0\) simultaneously evaluated at \([\bx_1,\dots,\bx_m] = [\bx_0+\boldsymbol{h}_1,\ldots,\bx_0+\boldsymbol{h}_m]\).
Indeed, let us put the simultaneous finite-difference evaluations into a matrix
\[
    \Delta_{\boldsymbol{H}} \bu(\bx_0) = [\bu(\bx_0+\boldsymbol{h}_1)-\bu(\bx_0),\ldots,\bu(\bx_0+\boldsymbol{h}_m)-\bu(\bx_0)] \in \R^{d\times m},
\]
then
\begin{equation}\label{eq:Taylor}
    \begin{aligned}
    \Delta_{\boldsymbol{H}} \bu(\bx_0) &= \boldsymbol{\mathcal{D}}\bu(\bx_0) \boldsymbol{H} + \boldsymbol{R}(\bx_0,\boldsymbol{H}),
    \quad\text{where}\\ 
    [\boldsymbol{\mathcal{D}}\bu]_{ij}(\bx_0) &= \partial_{x_j} u_i(\bx_0), \,i=1,\dots,d,\, j=1,\dots,N,\quad\text{and}\quad
    \lim_{\boldsymbol{H}\to\boldsymbol{0}} \boldsymbol{R}(\bx_0,\boldsymbol{H})|\boldsymbol{H}|_{F}^{-1} = \boldsymbol{0}, 
    \end{aligned}
\end{equation}
provided that \(\bu\) is continuously differentiable at \(\bx_0\).
We will consistently use \(|\cdot|_F\) to denote the Frobenius norm of the matrix.
Ignoring the higher-than-first order terms in \eqref{eq:Taylor}, we obtain a system of \(d\times m\) linear equations in terms of the \(d\times N\) unknown components of \(\boldsymbol{\mathcal{D}}\bu(\bx_0)\).
When \(m<N\) and \(\boldsymbol{H}\) is not a square matrix, it may be impossible to approximate the whole matrix \(\boldsymbol{\mathcal{D}}\bu(\bx_0)\).
Nevertheless, we can reconstruct the approximation of \(\boldsymbol{\mathcal{D}}\bu(\bx_0)\) restricted to the \(\operatorname{span}\{\boldsymbol{h}_1,\dots,\boldsymbol{h}_m\}=\operatorname{range}\boldsymbol{H}\) by computing
\begin{equation}\label{eq:def_simpl_gradient}
    \bD\boldsymbol{u}(\bx_0;\boldsymbol{H}) = \Delta_{\boldsymbol{H}}\bu(\bx_0)\boldsymbol{H}^\dagger \in \R^{d\times N},
\end{equation}
where \(\boldsymbol{H}^\dagger \in \R^{m \times N}\) is the Moore-Penrose inverse of \(\boldsymbol{H}\).
Indeed, representation~\eqref{eq:Taylor} shows that
\begin{equation}\label{eq:Taylor1}
    \bD\boldsymbol{u}(\bx_0;\boldsymbol{H}) = \boldsymbol{\mathcal{D}}\bu(\bx_0) \boldsymbol{H}\boldsymbol{H}^\dagger + \boldsymbol{R}(\bx_0,\boldsymbol{H})\boldsymbol{H}^\dagger.
\end{equation}
Since the matrix \(\boldsymbol{H}\boldsymbol{H}^\dagger\) is the orthogonal projection onto \(\operatorname{range}\boldsymbol{H}\), the first term in \eqref{eq:Taylor1} is the restriction of \(\boldsymbol{\mathcal{D}}\bu(\bx_0)\) to the subspace spanned by the simplex edges.
The remainder tends to zero provided that \(|\boldsymbol{H}|_F\longrightarrow 0\) while 
\(|\boldsymbol{H}|_F\,|\boldsymbol{H}^\dagger|_F\) stays bounded, that is, along uniformly
nondegenerate shrinking families of simplices.
By aggregating this quantity from a sufficiently broad family of matrices \(\boldsymbol{H}\) we may hope to obtain a reliable approximation 
to \(\boldsymbol{\mathcal{D}}\bu\).
In particular:
\begin{itemize}
\item If \(m=1\) we obtain a matrix containing the usual first order finite-difference quotient information, namely
\[
\begin{aligned}
    \bD\bu(\bx_0;[\boldsymbol{h}_1]) &= \frac{\bu(\bx_0+\boldsymbol{h}_1)-\bu(\bx_0)}{|\boldsymbol{h}_1|}\otimes \frac{\boldsymbol{h}_1}{|\boldsymbol{h}_1|},\quad\text{and therefore}\\
    |\bD\bu(\bx_0;[\boldsymbol{h}_1])|_F &= \frac{|\bu(\bx_0+\boldsymbol{h}_1)-\bu(\bx_0)|}{|\boldsymbol{h}_1|}, \quad\text{and, provided \(d=N\),}\\
    \operatorname{tr}\bD\bu(\bx_0;[\boldsymbol{h}_1]) &= \frac{\bu(\bx_0+\boldsymbol{h}_1)-\bu(\bx_0)}{|\boldsymbol{h}_1|}\cdot \frac{\boldsymbol{h}_1}{|\boldsymbol{h}_1|}.
\end{aligned}
\]
Integrating such quantities to obtain the energy of the system leads to well-known bond-based peridynamic models~\cite{Du2019, mengesha2015variational, Silling2000}.
\item If \(m=N\), \(\bu(\bx) = \boldsymbol{A}\bx + \boldsymbol{u}_0\) is an affine function, and \(\boldsymbol{H}=[\boldsymbol{h}_1,\dots,\boldsymbol{h}_N]\) is nonsingular (i.e., the simplex \(\operatorname{conv}\{ \bx, \bx+\boldsymbol{h}_1, \ldots, \bx+\boldsymbol{h}_{N} \}\) is nondegenerate) then \(\bD\boldsymbol{u}(\bx_0;\boldsymbol{H}) = \boldsymbol{A} = \boldsymbol{\mathcal{D}}\bu(\bx_0)\).
This is precisely the construction considered in~\cite{pedregal_special_2024}.
\end{itemize}

With this notation in place, we will consider nonlocal variational problems associated with minimizing energy functionals of the form
\begin{equation}\label{general}
    E(\bu) = \int_{\O} \int_{\mathcal{H}(\bx)} W(\bx,\boldsymbol{H},\bu(x),\bD\bu(\bx;\boldsymbol{H}))
    \,\mathrm{d}\boldsymbol{H}\,\mathrm{d}\bx,
\end{equation}
where \(\O \subset \R^N\) is a bounded regular domain,
\(W: \R^N\times \R^{N\times m} \times \R^d\times \R^{d\times N} \to \mathbb{R}\cup\{+\infty\}\) is an energy density,
\(\bu: \O \to \R^d\) are competing functions taken from a suitable family,
and
\[
    \mathcal{H}(\bx) = \{\, \boldsymbol{H}=[\boldsymbol{h}_1,\dots,\boldsymbol{h}_m] \in \R^{N\times m} \colon 
    \bx + \boldsymbol{h}_i \in \O, i=1,\dots,m \,\},
\]
is the set of admissible increment matrices \(\boldsymbol{H}\) keeping the finite difference calculations inside \(\O\).
For mechanical applications, it is often more natural to formulate the constitutive response in terms of the simplicial deformation gradient 
rather than the simplicial displacement gradient. 
If \(d=N\) and \[ \boldsymbol{y}(\bx)=\bx+\bu(\bx), \] we define 
\[ \boldsymbol{E}\bu( \bx;\boldsymbol{H} ) = \Delta_{\boldsymbol{H}}\boldsymbol{y}(\bx) \boldsymbol{H}^\dagger = \boldsymbol{H}\boldsymbol{H}^\dagger + \bD\bu(\bx;\boldsymbol{H}). \]
When \(m=N\) and \(\boldsymbol{H}\) is nonsingular, this reduces to 
\[ \boldsymbol{E}\bu( \bx;\boldsymbol{H} ) = \boldsymbol{I} + \bD\bu(\bx;\boldsymbol{H}), \] 
and represents the gradient of the unique affine deformation interpolating \(\boldsymbol{y}\) at the vertices of the interacting simplex.

This observation places continuum-kinematics-inspired peridynamics within a broader simplex-based constitutive framework.
In the original CPD formulation, separate one-, two-, and three-neighbour interaction energies are assigned to deformed line, area, and volume elements.
For a full-dimensional simplex, these quantities are generated, respectively, by the simplicial deformation gradient, its cofactor, and its determinant.
One may therefore replace the prescribed additive decomposition by a general stored-energy density assigned directly to the affine deformation of the interacting simplex.
We refer to the resulting class of models as \emph{simplicial hyperelasticity}.

The principal contributions of the present work are threefold. 

First, we show that bond-based peridynamics and the one-, two-, and three-neighbour interactions of continuum-kinematics-inspired peridynamics admit a common formulation in terms of simplicial deformation gradients. 
In the three-dimensional full-rank case, the associated line, area, and volume changes are described by the simplicial deformation gradient, its cofactor, and its determinant, respectively. 
This leads naturally to simplex-level hyperelastic energies with general, possibly nonadditive, constitutive dependence. 

Second, we introduce radial- and shape-weighted nonlocal spaces adapted to simplicial gradients. 
By passing from the multipoint interaction measure to its edge marginals, we obtain two-point estimates that yield nonlocal Poincar{\'e} inequalities and fractional Sobolev embeddings. 
For isotropic power-law interactions with at least two neighbours, the edge marginals exhibit a saturation effect. 

Third, we establish existence of equilibrium configurations through two complementary lower-semicontinuity mechanisms. 
Particle-wise averaged convexity yields weak lower semicontinuity without a compact embedding, whereas fractional compactness permits nonconvex energy densities.
The common mechanism behind the coercivity and compactness results is the reduction 
\[ \text{multipoint simplicial control} \quad\Longrightarrow\quad \text{two-point control through edge marginals}. \] 
The strength of the induced two-point interaction determines whether one obtains merely a nonlocal Poincar{\'e} inequality or the stronger fractional regularity needed for compactness.

The paper is organized as follows. 
In Section~\ref{existing-models}, we relate the general energy \eqref{general} to bond-based peridynamics and to the one-, two-, and three-neighbour interactions of CPD. 
Section~\ref{sec:simplicial-hyperelasticity} introduces simplicial hyperelasticity and discusses its principal constitutive requirements. 
In Section~\ref{sec:weighted-NW}, we define the radial- and shape-weighted nonlocal spaces, establish their basic functional-analytic properties, and develop the edge-marginal comparison, nonlocal Poincar{\'e} inequalities, and fractional Sobolev embeddings. 
Section~\ref{sec:existence} utilizes the results of Section~\ref{sec:weighted-NW} to develop the lower-semicontinuity and existence theory for the simplicial energy functionals we consider. 
Finally, Section~\ref{sec:measure-valued-perspective} places the Moore--Penrose construction within a broader parametrized-measure description of simplicial nonlocality and outlines directions for future work.

%=====================================================
% relation to existing models
\input{section_models.tex}
%=====================================================
% hyperelasticity
\input{section_hyperelasticity.tex}
%=====================================================
% weighted NW-spaces
\input{section_nwspaces.tex}
%=====================================================
% existence of solutions
\input{section_existence.tex}

\section{Concluding remarks and a measure-valued perspective} \label{sec:measure-valued-perspective}

We have built a variational and functional-analytic foundation for simplicial hyperelasticity and continuum-kinematics-inspired peridynamics. 
The line-, area-, and volume-interaction mechanisms of CPD have been expressed through a common simplicial deformation gradient; 
radial- and shape-weighted spaces have been introduced for the resulting energies; 
and the edge marginals of the interaction measure have been used to obtain coercivity and compactness. 
These results lead to two complementary existence theories, based respectively on particle-wise averaged convexity and fractional compactness. 
We conclude this work by placing the simplicial-gradient construction within a broader measure-valued perspective and indicating several directions for further research.

Inspired by the general framework of parametrized measures in nonlinear variational problems~\cite{Pedregal1997}, one can further generalize the
the simplicial-gradient construction presented here along the ideas outlined in~\cite{pedregal_special_2024}.
Let \( \boldsymbol{y}:\Omega\to\R^d \) be a measurable mapping; in the mechanical setting \(d=N\), we have \(\boldsymbol{y}(\bx)=\bx+\bu(\bx).\)
Consider an admissible increment matrix
\( \boldsymbol H = [\boldsymbol h_1,\ldots,\boldsymbol h_m]. \)
The associated multipoint data determine the action of an affine deformation map on the sampled interaction subspace through the constraint \(\boldsymbol{F}\boldsymbol{H} = \Delta_{\boldsymbol{H}}\boldsymbol{y}(\bx). \) 
Accordingly, define the affine manifold 
\begin{equation}
\mathcal{A}_{\bx,\boldsymbol{y}}(\boldsymbol{H}) = 
\left\{ \boldsymbol{F}\in\R^{d\times N}: \boldsymbol{F}\boldsymbol{H} = \Delta_{\boldsymbol{H}}\boldsymbol{y}(\bx) \right\}. 
\label{eq:compatible-affine-fibre} 
\end{equation} 
If \(\boldsymbol{H}\) has full column rank, every element of \(\mathcal{A}_{\bx,\boldsymbol{y}}(\boldsymbol{H})\) has the representation 
\begin{equation} \begin{aligned} 
\boldsymbol{F} &= \Delta_{\boldsymbol{H}}\boldsymbol{y}(\bx) \boldsymbol{H}^\dagger + \boldsymbol{Z} \bigl( \boldsymbol{I}_N - \boldsymbol{H}\boldsymbol{H}^\dagger \bigr), \qquad \boldsymbol{Z}\in\R^{d\times N}. 
\end{aligned} \label{eq:affine-fibre-parametrization} \end{equation} 
The Moore--Penrose construction employed throughout the manuscript selects the unique element of \(\mathcal{A}_{\bx,\boldsymbol{y}}(\boldsymbol{H})\) having minimal Frobenius norm, namely
\( \boldsymbol{F}_{\min}(\bx;\boldsymbol{H}) = \Delta_{\boldsymbol{H}}\boldsymbol{y}(\bx) \boldsymbol H^\dagger. \)
When \(m=N\) and \(\boldsymbol{H}\) is nonsingular, then \(\mathcal{A}_{\bx,\boldsymbol{y}}(\boldsymbol{H})\) is a singleton and this selection coincides with the complete simplicial 
deformation gradient. 
When \(m<N\), the second term in \eqref{eq:affine-fibre-parametrization} describes the undetermined action on the orthogonal complement of \(\operatorname{range}\boldsymbol{H}\). 
A measure-valued extension of the model may be obtained by assigning to almost every admissible pair \((\bx,\boldsymbol{H})\) a probability measure 
\( \xi_{\bx,\boldsymbol{y},\boldsymbol{H}} \in \mathcal{P}(\R^{d\times N}) \)
such that 
\begin{equation} 
\operatorname{supp} \xi_{\bx,\boldsymbol{y},\boldsymbol{H}} \subset \mathcal{A}_{\bx,\boldsymbol{y}}(\boldsymbol{H}). 
\label{eq:kernel-support-affine-fibre} 
\end{equation} 
Then
\( (\bx,\boldsymbol{H}) \longmapsto \xi_{\bx,\boldsymbol{y},\boldsymbol{H}} \) 
describes the distribution of affine deformation maps compatible with the multipoint observations associated with \(\boldsymbol{H}\);
we assume that this kernel is  weakly measurable.
Let \(\Lambda_{\bx}\) denote the interaction measure on admissible increment matrices induced by the radial and shape-orientation measures. 
If \( 0< \Lambda_{\bx}\bigl(\mathcal{H}(\bx)\bigr) <+\infty, \) we can normalize this measure by 
\[ \widehat{\Lambda}_{\bx} = \frac{ \Lambda_{\bx} }{ \Lambda_{\bx}(\mathcal{H}(\bx)) }. \] 
If the total interaction measure is infinite, the same construction may be applied to finite truncations, or formulated using nonnegative parametrized measures without normalization.
Averaging the conditional measures over the interaction configurations yields 
\begin{equation} \begin{aligned} 
\lambda_{\bx,\boldsymbol{y}}(B) = \int_{\mathcal H(\bx)} \xi_{\bx,\boldsymbol{y},\boldsymbol{H}}(B) \,\mathrm{d}\widehat{\Lambda}_{\bx}(\boldsymbol{H}),
\end{aligned} \label{eq:parametrized-affine-measure} \end{equation} 
for every Borel set \(B\subset\R^{d\times N}\). 
The family 
\( \bigl( \lambda_{\bx,\boldsymbol{y}} \bigr)_{\bx\in\Omega} \)
is a parametrized probability measure of Young-measure type. 
In contrast to a classical Young measure generated by oscillations of a sequence, it describes the distribution of affine deformation maps induced by 
the family of nonlocal interaction configurations at a fixed deformation \(\boldsymbol{y}\). 
The model studied in this work corresponds to the Dirac measure setting 
\begin{equation} \xi_{\bx,\boldsymbol{y},\boldsymbol{H}} = \delta_{ \Delta_{\boldsymbol{H}}\boldsymbol{y}(\bx) \boldsymbol H^\dagger }, \label{eq:Moore-Penrose-Dirac-selection} \end{equation} 
where
\[ \lambda_{\bx,\boldsymbol{y}} = \left( \boldsymbol{H} \longmapsto \Delta_{\boldsymbol{H}}\boldsymbol{y}(\bx) \boldsymbol H^\dagger \right)_\# \widehat\Lambda_{\bx}. \] 
Thus, the simplicial-gradient formulation amounts to concentrating the measure on the minimum-norm affine map compatible with each interaction configuration. 
If the stored-energy density depends explicitly on \(\boldsymbol{H}\), the marginal measure \(\lambda_{\bx,\boldsymbol{y}}\) does not retain all constitutively relevant information. 
In that case, one may retain the joint measure 
\begin{equation} \begin{aligned} 
\gamma_{\bx,\boldsymbol{y}} ( \mathrm{d}\boldsymbol{H}, \mathrm{d}\boldsymbol{F} ) = \xi_{\bx,\boldsymbol{y},\boldsymbol{H}} (\mathrm{d}\boldsymbol{F}) \widehat\Lambda_{\bx} (\mathrm{d}\boldsymbol{H}). 
\end{aligned} \label{eq:joint-parametrized-measure} \end{equation} 
The corresponding local contribution to the energy takes the form 
\[ 
\int W( \bx, \boldsymbol{H}, \boldsymbol{y}(\bx), \boldsymbol{F} ) \,\mathrm{d}\gamma_{\bx,\boldsymbol{y}} ( \boldsymbol{H}, \boldsymbol{F} ). 
\] 
The measure-valued extension separates two ingredients of the nonlocal description. 
The interaction measure specifies the distribution of scales, orientations, and simplex geometries, whereas the conditional kernel describes 
the affine deformation maps compatible with each sampled configuration. 
Developing this larger model would require suitable measurability and compatibility conditions, a characterization of the parametrized measures 
generated by deformations, and a corresponding relaxation and lower-semicontinuity theory. 

Other natural questions concern the dependence of the resulting spaces and constitutive classes on the interaction order \(m\), the localization limits of simplicial energies, and the characterization of classical elastic and hyperelastic responses obtainable from different families of interaction simplices. Discrete shape-orientation measures also provide a direct connection with numerical models based on finite collections of simplex configurations. These directions are left for future work.

\printbibliography
\end{document}

%% file: section_models.tex
%=====================================================
\section{Relation to existing nonlocal models of mechanics}
\label{existing-models}

We will now relate the model~\eqref{general} to existing nonlocal models of continuum mechanics, showing that it indeed extends typical nonlocal functionals from bond-based peridynamics as well as the two- and three-neighbour interactions of CPD models.

\subsection{Bond-based peridynamic diffusion}
\label{subsec:bb_pd_diffusion}

Scalar peridynamic diffusion may be easily cast into the form~\eqref{general} as follows.
We begin by setting \(d=m=1\).
Let \(\delta>0\) be a nonlocal horizon, determining the largest distance at which two points may interact, and let
\(\hat{\rho}: [0,\infty)\to [0,\infty)\) be a radial kernel, determining the strength of such interactions.
Thus one assumes \(\hat{\rho}(r) = 0\), for \(r>\delta\), and often (but not always) a normalization condition
\(
    \int_{\R^N}\hat{\rho}(|\boldsymbol{h}|)\,\mathrm{d}\boldsymbol{h} = 1.
\)
We further let
\[
 W(\bx,\boldsymbol{H},u,\boldsymbol{G}) = \hat{\rho}(|\boldsymbol{H}|_F)W_1(|\boldsymbol{G}|_F),
\]
where \(W_1: \R\to [0,\infty)\) is a pairwise bond-interaction energy density.
For example, in the paradigmatic nonlocal approximation of \(p\)-Laplacian we have \(W_1(t) = p^{-1} |t|^p\) for some \(p\in (1,\infty)\).
After a change of integration variables from \((\bx,\boldsymbol{H})\) to \((\bx_0,\bx_1)\) in~\eqref{general} we obtain a familiar form of a bond-based peridynamic energy
\[
    E(u) = \int_{\O}\int_{\O} \hat{\rho}(|\bx_1-\bx_0|)W_1\bigg(\frac{|u(\bx_1)-u(\bx_0)|}{|\bx_1-\bx_0|}\bigg)\,\mathrm{d}\bx_1\,\mathrm{d}\bx_0.
\]

\subsection{Bond-based peridynamic solid}
\label{subsec:bb_pd_solid}

This case is essentially the same as that of the scalar diffusion described in Subsection~\ref{subsec:bb_pd_diffusion},
except instead of a scalar-valued function corresponding to \(d=1\) we have a vector-valued displacement function corresponding to \(d=N\).
Having introduced the radial kernel \(\hat{\rho}\) and a pairwise bond-interaction energy density \(W_1\) as above, we put
\[
 W(\bx,\boldsymbol{H},u,\boldsymbol{G}) = \hat{\rho}(|H|_F)W_1(\operatorname{tr} \boldsymbol{G}).
\]
This results in the familiar model~\cite{mengesha2015variational}:
\begin{equation}\label{two-point}
    E(\bu) = \int_{\O}\int_{\O}
    \hat{\rho}(|\bx_1-\bx_0|)
    W_1\bigg(
    \frac{\bu(\bx_1)-\bu(\bx_0)}{|\bx_1-\bx_0|}
    \cdot \frac{\bx_1-\bx_0}{|\bx_1-\bx_0|}
    \bigg)\,\mathrm{d}\bx_1\,\mathrm{d}\bx_0.
\end{equation}

\subsection{Continuum-kinematics-inspired peridynamics (CPD)}
\label{subsec:cpd} 

We next relate the simplicial gradient to continuum-kinematics-inspired peridynamics (CPD), introduced in~\cite{javili2019continuum}.
For simplicity, we restrict the discussion in this subsection to the mechanically most relevant case \(d=N=m=3.\)
Thus, \(\bu:\O\to\R^3\) is a displacement, and  \(\boldsymbol{y}(\bx)=\bx+\bu(\bx)\) is the corresponding deformation.

\paragraph{Line, area, and volume elements in CPD.} 
Let us fix \(\bx\in\O\), and consider three neighbouring material points \(\bx_i = \bx+\boldsymbol{h}_i\), \(i=1,2,3,\) such that 
\(\boldsymbol{H} = [\boldsymbol{h}_1,\boldsymbol{h}_2,\boldsymbol{h}_3]\in\R^{3\times 3} \) is nonsingular. 
Following the terminology of~\cite{javili2019continuum}, the vectors \(\boldsymbol{h}_i\) represent finite relative-position vectors in the reference configuration, while the vectors
\begin{equation} \label{eq:cpd-deformed-bonds}
    \boldsymbol{\xi}_i = \boldsymbol{y}(\bx+\boldsymbol{h}_i)-\boldsymbol{y}(\boldsymbol{x}) = 
    \boldsymbol{h}_i+\bu(\bx+\boldsymbol{h}_i)-\bu(\bx), \qquad i=1,2,3,
\end{equation} 
represent their counterparts in the deformed configuration.
The elementary CPD kinematic quantities are constructed from these relative vectors as follows: 
\begin{equation} 
    \label{eq:cpd-kinematic-quantities}
    \begin{aligned} 
        &\text{one-neighbour quantity:} &&\boldsymbol{\xi}_i, \\[1mm] 
        &\text{two-neighbour quantity:} &&\boldsymbol{a}_{ij}=\boldsymbol{\xi}_i\times\boldsymbol{\xi}_j, \\[1mm] 
        &\text{three-neighbour quantity:} &&v_{ijk}=(\boldsymbol{\xi}_i\times\boldsymbol{\xi}_j)\cdot\boldsymbol{\xi}_k. 
    \end{aligned} 
\end{equation}
The vector \(\boldsymbol{a}_{ij}\) is the oriented area vector of the parallelogram spanned by the two deformed bonds, while \(v_{ijk}\) 
is the oriented volume of the parallelepiped spanned by the three deformed bonds.\footnote{The corresponding triangle and tetrahedron have oriented area vector \(\boldsymbol{a}_{ij}/2\) and oriented volume \(v_{ijk}/6\), respectively. The factors \(1/2\) and \(1/6\) play no
essential role in the following kinematic identities and can be incorporated into the interaction densities.}
The reference configuration counterparts of these quantities are 
\begin{equation} \label{eq:cpd-reference-area-volume} 
    \boldsymbol{h}_i,\qquad
    \boldsymbol{A}_{ij}=\boldsymbol{h}_i\times \boldsymbol{h}_j, 
    \quad\text{and}\quad V_{ijk}=(\boldsymbol{h}_i\times \boldsymbol{h}_j)\cdot \boldsymbol{h}_k. 
\end{equation}

\paragraph{The original CPD energy.} 
The internal energy in the original CPD formulation~\cite{javili2019continuum} makes a
careful distinction between one-, two-, and three-neighbour interaction energies.
Specifically, it is postulated as an additive combination 
\(E_{\text{CPD}}(\bu) = E_1(\bu) + E_2(\bu) + E_3(\bu)\),
where the three terms describe interactions between a material point and, respectively, 
one, two, or three of its neighbours. 
These take the form of double, triple, and quadruple integrals. 
In the notation introduced above, the three energy contributions have the conceptual form\footnote{The factors \(1/2\), \(1/3\), and \(1/4\) follow the convention of~\cite{javili2019continuum}; equivalently, they may be absorbed into the interaction densities.}
\[
    \begin{aligned}
    E_1(\bu) &= \frac{1}{2}\int_\O \int_{\O-\bx} w_1\bigl( \boldsymbol{\xi}_1 \bigr) \,\mathrm{d}\boldsymbol{h}_1\,\mathrm{d}\boldsymbol{x},\\
    E_2(\bu) &= \frac{1}{3} \int_\O \int_{\O-\bx} \int_{\O-\bx} w_2\bigl( \boldsymbol{\xi}_1\times\boldsymbol{\xi}_2 \bigr) \,\mathrm{d}\boldsymbol{h}_2\,\mathrm{d}\boldsymbol{h}_1\,\mathrm{d}\bx, \quad\text{and}\\
    E_3(\bu) &= \frac{1}{4} \int_\O \int_{\O-\bx} \int_{\O-\bx} \int_{\O-\bx} w_3\bigl( (\boldsymbol{\xi}_1\times\boldsymbol{\xi}_2)\cdot\boldsymbol{\xi}_3 \bigr) \,\mathrm{d}\boldsymbol{h}_3\,\mathrm{d}\boldsymbol{h}_2\,\mathrm{d}\boldsymbol{h}_1\,\mathrm{d}\boldsymbol{x},
    \end{aligned}
\]
where the energy densities depend on \(\bu\) and the integration variables via~\eqref{eq:cpd-deformed-bonds} and~\eqref{eq:cpd-kinematic-quantities}.
As in Subsection~\ref{subsec:bb_pd_diffusion} one may also include kernels accounting for
varying interaction strengths and a finite interaction horizon.

\paragraph{The simplicial deformation gradient.} 
To fit the CPD energy into our unified framework, we proceed as follows.
Let us collect the deformed relative vectors into the matrix \(\Delta_{\boldsymbol{H}}\boldsymbol{y}(\bx) = [\boldsymbol{\xi}_1,\boldsymbol{\xi}_2,\boldsymbol{\xi}_3].\) 
Owing to~\eqref{eq:cpd-deformed-bonds} we can express it as \(\Delta_{\boldsymbol{H}}\boldsymbol{y}(\bx) = \boldsymbol{H}+\Delta_{\boldsymbol{H}}\bu(\bx).\)
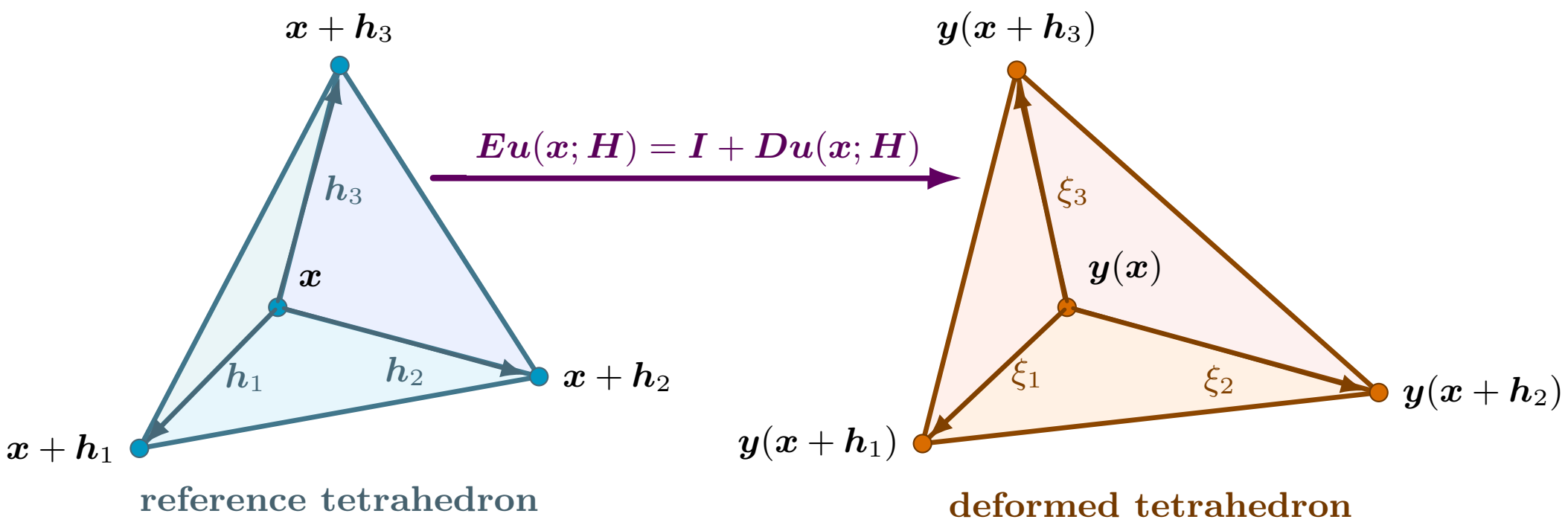
\begin{figure}[t] \centering 
    \input{simplex_tikz_figure}
\caption{Reference and deformed tetrahedra associated with the simplicial deformation gradient. The deformation of line, oriented area, and oriented volume elements is governed by \(\boldsymbol{E}\bu(\bx;\boldsymbol{H})\), \(\operatorname{cof} \boldsymbol{E}\bu(\bx;\boldsymbol{H})\), and \(\det \boldsymbol{E}\bu(\bx;\boldsymbol{H})\), respectively; see~\eqref{eq:cpd-line-map}, \eqref{eq:cpd-area-cofactor}, and~\eqref{eq:cpd-volume-determinant}.} 
\label{fig:simplicial-deformation} 
\end{figure}
The unique affine deformation interpolating the values of \(\boldsymbol{y}\) at the four vertices of the tetrahedron has gradient 
\begin{equation}  \label{eq:cpd-simplicial-deformation-gradient}
  \begin{aligned} 
    \boldsymbol{E}\bu(\bx;\boldsymbol{H}) = \Delta_{\boldsymbol{H}}\boldsymbol{y}(\bx)\boldsymbol{H}^{-1}
    = \bigl(\boldsymbol{H}+\Delta_{\boldsymbol{H}}\bu(\bx)\bigr)\boldsymbol{H}^{-1} = 
    \boldsymbol{I}+\boldsymbol{D}\bu(\bx;\boldsymbol{H}), 
  \end{aligned}
\end{equation} 
as \(\boldsymbol{H}^\dagger=\boldsymbol{H}^{-1}\) in the present full-rank setting.
We can refer to \(\boldsymbol{E}\bu(\bx;\boldsymbol{H})\) as the \emph{simplicial deformation gradient}, see Fig.~\ref{fig:simplicial-deformation}.
In particular, we have the following identity relating one-neighbour quantities
\begin{equation}\label{eq:cpd-line-map}
\boldsymbol{E}\bu(\bx;\boldsymbol{H})\boldsymbol{h}_i=\boldsymbol{\xi}_i, \qquad i=1,2,3.
\end{equation}
The distinction between \(\boldsymbol{D}\bu(\cdot;\boldsymbol{H})\) and \(\boldsymbol{E}\bu(\cdot,\boldsymbol{H})\) is important insofar as 
the former describes displacement differences, whereas the latter describes the deformation of the material tetrahedron.
For example, when \(\bu=\boldsymbol{0}\), one has \(\boldsymbol{D}\bu(\cdot;\boldsymbol{H})=\boldsymbol{0}\) but \(\boldsymbol{E}\bu(\cdot;\boldsymbol{H})=\boldsymbol{I}\), as the reference line, area, and volume elements remain unchanged.
If \(\bu\) is affine, say \(\bu(x)=\boldsymbol{G}\bx+\bu_0, \) then \(\boldsymbol{D}\bu(\bx;\boldsymbol{H})=\boldsymbol{G}\),
while \( \boldsymbol{E}\bu(\bx;\boldsymbol{H})=\boldsymbol{I}+\boldsymbol{G}\), for every \(\bx\) and every non-singular \(\boldsymbol{H}\).
Thus in the affine case the simplicial deformation gradient agrees exactly with the classical deformation gradient \(\boldsymbol{\mathcal{D}}\boldsymbol{y}=\boldsymbol{I}+\boldsymbol{\mathcal{D}}\bu.\)

\paragraph{Two-neighbour kinematics and the cofactor.}
For every \(\boldsymbol{F}\in\R^{3\times3}\) and every \(\boldsymbol{a},\boldsymbol{b}\in\R^{3}\), the cross-product identity 
\(
    (\boldsymbol{F}\boldsymbol{a})\times(\boldsymbol{F}\boldsymbol{b}) = (\operatorname{cof} \boldsymbol{F})(\boldsymbol{a}\times \boldsymbol{b})  
\) 
holds. 
Applying this identity to \(\boldsymbol{F}=\boldsymbol{E}\bu(\bx;\boldsymbol{H})\), \(\boldsymbol{a}=\boldsymbol{h}_i\), \(\boldsymbol{b}=\boldsymbol{h}_j\), and using \eqref{eq:cpd-line-map}, we obtain the relation 
\begin{equation} \label{eq:cpd-area-cofactor}
    \boldsymbol{a}_{ij} = \boldsymbol{\xi}_i\times\boldsymbol{\xi}_j 
    = 
    \bigl(\boldsymbol{E}\bu(\bx;\boldsymbol{H})\boldsymbol{h}_i\bigr) \times \bigl(\boldsymbol{E}\bu(\bx;\boldsymbol{H})\boldsymbol{h}_j\bigr)
    = 
    \operatorname{cof}\bigl(\boldsymbol{E}\bu(\bx;\boldsymbol{H})\bigr)\boldsymbol{A}_{ij}. 
\end{equation} 
Consequently, the CPD two-neighbour quantity is the image of the reference oriented area vector under the cofactor of the simplicial 
deformation gradient. 
This is the finite-scale counterpart of Nanson's formula in classical continuum mechanics;\footnote{%
Note that a pair of bonds \(\boldsymbol{h}_i,\boldsymbol{h}_j\) determines the reference triangle but does not determine the action of a linear map on the whole of \(\R^3\). 
Nevertheless, any linear map \(\boldsymbol{F}\) satisfying 
\(\boldsymbol{F}\boldsymbol{h}_i=\boldsymbol{\xi}_i\), \(\boldsymbol{F}\boldsymbol{h}_j=\boldsymbol{\xi}_j \)
also necessarily satisfies 
\((\operatorname{cof} \boldsymbol{F})(\boldsymbol{h}_i\times \boldsymbol{h}_j)=\boldsymbol{\xi}_i\times\boldsymbol{\xi}_j. \)}
see, for example, \cite[Sections~1.5--1.8]{Ciarlet1988} or \cite[Chapter~2]{Holzapfel2000}.

\paragraph{Three-neighbour kinematics and the determinant.} 
For every \(\boldsymbol{F}\in\R^{3\times3}\) and every \(\boldsymbol{a},\boldsymbol{b},\boldsymbol{c}\in\R^3\), one has the scalar
triple-product identity 
\(
  \bigl((\boldsymbol{F}\boldsymbol{a})\times(\boldsymbol{F}\boldsymbol{b})\bigr)\cdot(\boldsymbol{F}\boldsymbol{c}) = 
  (\det \boldsymbol{F})(\boldsymbol{a}\times \boldsymbol{b})\cdot \boldsymbol{c}. 
  \label{eq:cpd-triple-product-identity} 
\) 
Applying this identity to the three reference bonds yields 
\begin{equation}
    v_{ijk} = (\boldsymbol{\xi}_i\times\boldsymbol{\xi}_j)\cdot\boldsymbol{\xi}_k 
    = \det\bigl(\boldsymbol{E}\bu(\bx;\boldsymbol{H})\bigr) V_{ijk}.
  \label{eq:cpd-volume-determinant} 
\end{equation}
In particular, \(\det \boldsymbol{E}\bu(\bx;\boldsymbol{H})=v_{123}/V_{123}\) is the signed ratio of the deformed and reference tetrahedral volumes.

\paragraph{The original CPD energy in the unified framework~\eqref{general}.}

With the help of the notation introduced above, we can express one-, two-, and three-neighbour interaction energies as follows:
\[
\begin{aligned}
E_1(\bu) &= \frac{1}{2|\Omega|^2}\int_\O \int_{\mathcal{H}(\bx)} w_1\bigl(  [\boldsymbol{I}+\boldsymbol{D}\bu(\bx;\boldsymbol{H})]\boldsymbol{H}\boldsymbol{e}_1 \bigr) \,\mathrm{d}\boldsymbol{H}\,\mathrm{d}\bx,\\
E_2(\bu) &= \frac{1}{3|\Omega|}\int_\O \int_{\mathcal{H}(\bx)} w_2\bigl( \operatorname{cof}\bigl(\boldsymbol{I}+\boldsymbol{D}\bu(\bx;\boldsymbol{H})\bigr)[(\boldsymbol{H}\boldsymbol{e}_1) \times (\boldsymbol{H}\boldsymbol{e}_2)] \bigr) \,\mathrm{d}\boldsymbol{H}\,\mathrm{d}\bx, \quad\text{and}\\
E_3(\bu) &= \frac{1}{4}\int_\O \int_{\mathcal{H}(\bx)}  w_3\bigl( \det\bigl(\boldsymbol{I}+\boldsymbol{D}\bu(\bx;\boldsymbol{H})\bigr) [(\boldsymbol{H}\boldsymbol{e}_1) \times (\boldsymbol{H}\boldsymbol{e}_2)]\cdot \boldsymbol{H}\boldsymbol{e}_3 \bigr) \,\mathrm{d}\boldsymbol{H}\,\mathrm{d}\bx,
\end{aligned}
\]
where \(\boldsymbol{e}_1,\boldsymbol{e}_2,\boldsymbol{e}_3\) denote the canonical basis vectors in \(\R^3\), so that \(\boldsymbol{H}\boldsymbol{e}_i=\boldsymbol{h}_i\).
Consequently, the original CPD model is contained in the general framework~\eqref{general}.
In the present three-dimensional setting, all three interaction terms may be represented 
on the common configuration space of tetrahedra.
The one- and two-neighbour terms are lifted to this space by integration with respect to 
the unused vertices, producing the normalization factors \(\lvert\O\rvert^{-2}\) and 
\(\lvert\O\rvert^{-1}\), respectively, while the three-neighbour term already has the
required quadruple-integral form. 
These normalization factors apply to unrestricted integration over \((\Omega-\bx)^3\). If localization kernels are included, the corresponding marginal weights must instead be incorporated into the lifted interaction densities.
The preceding representation also reveals that the deformed line, area, and volume elements
used in CPD are generated by the simplicial deformation gradient 
\(\boldsymbol{E}\bu(\bx;\boldsymbol{H})\), its cofactor, and its determinant.
This observation suggests a broader constitutive framework in which the energy is assigned
directly to the affine deformation of each interacting simplex. 
We develop this proposal in the following section. 
Even more generally, one may replace the single affine deformation map selected by the Moore--Penrose inverse with a parametrized probability measure on the family of affine maps compatible with the multipoint data. 
This yields a measure-valued extension of the present framework, briefly discussed in Section~\ref{sec:measure-valued-perspective}.

%% file: simplex_tikz_figure.tex
    \tdplotsetmaincoords{68}{112} 
    \begin{tikzpicture}[ tdplot_main_coords, scale=1.15, line cap=round, line join=round, >=Latex, refedge/.style={ draw=cyan!45!black, line width=1.25pt }, defedge/.style={ draw=orange!55!black, line width=1.25pt }, hidden/.style={ dashed, opacity=.65, line width=.8pt }, refpoint/.style={ circle, draw=cyan!35!black, fill=cyan!75!black, inner sep=1.6pt }, defpoint/.style={ circle, draw=orange!45!black, fill=orange!85!black, inner sep=1.6pt }, maparrow/.style={ -{Latex[length=3.2mm,width=2.1mm]}, draw=violet!75!black, line width=1.6pt }, every node/.style={font=\small} ]
    % ------------------------------------------------- 
    % Reference tetrahedron 
    % ------------------------------------------------- 
    \begin{scope}
      \coordinate (X0) at (0,0,0); 
      \coordinate (X1) at (3.1,0.1,0.0); 
      \coordinate (X2) at (0.8,2.5,0.1); 
      \coordinate (X3) at (0.7,0.8,2.4); 
      % translucent faces, back to front 
      \fill[cyan!18,opacity=.42] (X0)--(X1)--(X2)--cycle; 
      \fill[blue!18,opacity=.34] (X0)--(X2)--(X3)--cycle; 
      \fill[teal!20,opacity=.34] (X0)--(X3)--(X1)--cycle; 
      \fill[cyan!10,opacity=.25] (X1)--(X2)--(X3)--cycle; 
      % edges 
      \draw[refedge] (X0)--(X1); 
      \draw[refedge] (X0)--(X2); 
      \draw[refedge] (X0)--(X3); 
      \draw[refedge] (X1)--(X2); 
      \draw[refedge] (X1)--(X3); 
      \draw[refedge] (X2)--(X3); 
      % points and labels 
      \node[refpoint,label=above right:$\bx$] at (X0) {}; 
      \node[refpoint,label=left:$\bx+\boldsymbol{h}_1$] at (X1) {}; 
      \node[refpoint,label=right:$\bx+\boldsymbol{h}_2$] at (X2) {}; 
      \node[refpoint,label=above:$\bx+\boldsymbol{h}_3$] at (X3) {}; 
      % edge vectors 
      \draw[-{Latex[length=2.3mm]},cyan!35!black,line width=1pt] ($(X0)+(0.05,0,0)$) -- node[pos=.5,right] {$\boldsymbol{h}_1$} ($(X1)+(-0.10,0,0)$); 
      \draw[-{Latex[length=2.3mm]},cyan!35!black,line width=1pt] ($(X0)+(0,0.05,0)$) -- node[pos=.5,below] {$\boldsymbol{h}_2$} ($(X2)+(0,-0.10,0)$); 
      \draw[-{Latex[length=2.3mm]},cyan!35!black,line width=1pt] ($(X0)+(0,0,0.05)$) -- node[pos=.5,right] {$\boldsymbol{h}_3$} ($(X3)+(0,0,-0.10)$); 
      \node[ anchor=north, text=cyan!30!black, font=\small\bfseries ] at ($(X1)!0.5!(X2)+(0,0,-0.5)$) {reference tetrahedron};
    \end{scope}
    % ------------------------------------------------- 
    % Deformed tetrahedron, translated to the right 
    % ------------------------------------------------- 
    \begin{scope}[xshift=6.1cm] 
        \coordinate (Y0) at (0,0,0); 
        \coordinate (Y1) at (3.6,0.25,0.25); 
        \coordinate (Y2) at (0.5,2.8,-0.10); 
        \coordinate (Y3) at (2.15,0.45,2.85); 
        % translucent faces 
        \fill[orange!22,opacity=.44] (Y0)--(Y1)--(Y2)--cycle; 
        \fill[magenta!18,opacity=.34] (Y0)--(Y2)--(Y3)--cycle; 
        \fill[red!18,opacity=.32] (Y0)--(Y3)--(Y1)--cycle; 
        \fill[orange!12,opacity=.26] (Y1)--(Y2)--(Y3)--cycle; 
        % edges 
        \draw[defedge] (Y0)--(Y1); 
        \draw[defedge] (Y0)--(Y2); 
        \draw[defedge] (Y0)--(Y3); 
        \draw[defedge] (Y1)--(Y2); 
        \draw[defedge] (Y1)--(Y3); 
        \draw[defedge] (Y2)--(Y3); 
        % points and labels 
        \node[defpoint,label=above right:$\boldsymbol{y}(\bx)$] at (Y0) {}; 
        \node[defpoint,label=left:$\boldsymbol{y}(\bx+\boldsymbol{h}_1)$] at (Y1) {}; 
        \node[defpoint,label=right:$\boldsymbol{y}(\bx+\boldsymbol{h}_2)$] at (Y2) {}; 
        \node[defpoint,label=above:$\boldsymbol{y}(\bx+\boldsymbol{h}_3)$] at (Y3) {}; 
        % deformed edge vectors 
        \draw[-{Latex[length=2.3mm]},orange!50!black,line width=1pt] ($(Y0)+(0.05,0,0)$) -- node[pos=.5,right] {$\xi_1$} ($(Y1)+(-0.10,0,0)$); 
        \draw[-{Latex[length=2.3mm]},orange!50!black,line width=1pt] ($(Y0)+(0,0.05,0)$) -- node[pos=.5,below] {$\xi_2$} ($(Y2)+(0,-0.10,0)$); 
        \draw[-{Latex[length=2.3mm]},orange!50!black,line width=1pt] ($(Y0)+(0,0,0.05)$) -- node[pos=.5,right] {$\xi_3$} ($(Y3)+(0,0,-0.10)$); 
        \node[ anchor=north, text=orange!45!black, font=\small\bfseries ] at ($(Y1)!0.5!(Y2)+(0,0,-0.5)$) {deformed tetrahedron}; 
    \end{scope} 
    % ------------------------------------------------- 
    % Central deformation arrow, in screen coordinates 
    % ------------------------------------------------- 
    \begin{scope}[tdplot_screen_coords] 
         \draw[maparrow] (1.2,1.0) -- node[ above, align=center, text=violet!70!black, fill=white, inner sep=2pt ] {$\boldsymbol{E}\bu(\bx;\boldsymbol{H})=\boldsymbol{I}+\boldsymbol{D}\bu(\bx;\boldsymbol{H})$} (5.3,1.0); 
        %\node[ draw=violet!30, fill=violet!4, rounded corners=2pt, inner sep=5pt, align=left, font=\small ] at (3.0,-2.5) { $\xi_i=F_u h_i$\\[1.5mm] $\xi_i\times\xi_j =(\operatorname{cof} F_u)(h_i\times h_j)$\\[1.5mm] $(\xi_1\times\xi_2)\cdot\xi_3 =(\det F_u)(h_1\times h_2)\cdot h_3$ }; 
    \end{scope}
\end{tikzpicture} 

%% file: section_hyperelasticity.tex
\section{Simplicial hyperelasticity} 
\label{sec:simplicial-hyperelasticity}

Much of the following discussion extends to \(N=d\geq2\) and \(1\leq m\leq N\). For clarity, however, we continue to focus on the mechanically most relevant case \(N=d=m=3\).
Recall that the complete energy density for the original CPD model~\cite{javili2019continuum} assigns separate interaction energies 
to deformed line, area, and volume elements. 
In this section we consider a significantly more general class of energies.

\paragraph{Simplicial stored energy densities.}
The simplicial deformation gradient permits a different construction: one may associate 
a stored-energy density with the deformation associated with the simplex. 
One may therefore consider energies of the form
\begin{equation} 
    E(\bu) = \int_\O \int_{\mathcal{H}(\bx)} \Psi\bigl(\bx,\boldsymbol{H},\boldsymbol{E}\bu(\bx;\boldsymbol{H})\bigr) \,\mathrm{d}\boldsymbol{H}\,\mathrm{d}\bx.\label{eq:simplicial-hyperelastic-energy} 
\end{equation}
To obtain the appropriate physical dimensions, a classical stored energy per unit reference volume may be multiplied by an interaction weight. For example, one may take 
\( \Psi(\bx,\boldsymbol{H},\boldsymbol{F}) = \rho(\bx,\boldsymbol{H})\Psi_{\rm loc}(\boldsymbol{F}), \)
where the normalization and physical units of \(\rho\) depend on the chosen family of interacting simplices.
This is a particular case of~\eqref{general}, obtained by setting 
\(W(\bx,\boldsymbol{H},\boldsymbol{z},\boldsymbol{G}) =\Psi(\bx,\boldsymbol{H},\boldsymbol{I}+\boldsymbol{G}).\)
Unlike the prescribed additive decomposition of the original CPD energy, 
\eqref{eq:simplicial-hyperelastic-energy} permits a general coupled dependence on the 
affine deformation of the simplex. 
The dependence on \(\bx\) and \(\boldsymbol{H}\) may describe material inhomogeneity, the scale and shape of the reference simplex, localization of the interactions, and material anisotropy.
The additive decomposition of the original CPD model is recovered only when \(\Psi\) can 
be expressed as a sum of contributions associated separately with selected deformed line,
area, and volume elements.

\paragraph{Frame indifference.}
In analogy with local hyperelasticity, a natural constitutive requirement is material frame indifference, which requires that the stored energy is unaffected by a superposed rigid motion of the deformed configuration, while the reference configuration remains fixed.
Here we restrict attention to the case $N=m$.
Thus let the deformed configuration transform as
\[ \boldsymbol{y}^*(\bx)=\boldsymbol{c}+\boldsymbol{Q}\boldsymbol{y}(\bx), 
\qquad \boldsymbol{Q}\in\mathrm{SO}(N), \boldsymbol{c}\in\R^N.
\] 
Since the reference configuration remains fixed, both \(\bx\) and \(\boldsymbol H\) remain unchanged.
The corresponding displacement is 
\[ \bu^*(\bx) = \boldsymbol{y}^*(\bx)-\bx 
= \boldsymbol{c}+(\boldsymbol{Q}-\boldsymbol{I})\bx+\boldsymbol{Q}\bu(\bx). 
\] 
It follows that 
\[ 
\boldsymbol{I}+\boldsymbol{D}\bu^*(\bx;\boldsymbol{H}) = 
\boldsymbol{Q}\bigl(\boldsymbol{I}+\boldsymbol{D}\bu(\bx;\boldsymbol{H})\bigr),
\quad\text{or equivalently,}\quad
    \boldsymbol{E}\bu^*(\bx;\boldsymbol{H})=
    \boldsymbol{Q}\boldsymbol{E}\bu(\bx;\boldsymbol{H}). 
\] 
Consequently, material frame indifference requires 
\begin{equation} 
\Psi(\bx,\boldsymbol{H},\boldsymbol{Q}\boldsymbol{F})=\Psi(\bx,\boldsymbol{H},\boldsymbol{F}),
    \quad\boldsymbol{Q}\in\mathrm{SO}(N). 
    \label{eq:simplicial-frame-indifference}
\end{equation}
As in classical nonlinear elasticity, on the orientation-preserving class 
\[ \mathrm{GL}^+(N) = \left\{ \boldsymbol{F}\in\R^{N\times N}: \det\boldsymbol{F}>0 \right\}, \]
condition~\eqref{eq:simplicial-frame-indifference} allows the deformation dependence to 
be expressed through the simplicial right Cauchy--Green tensor 
\begin{equation} 
  \boldsymbol{C}\bu(\bx;\boldsymbol{H}) = \boldsymbol{E}\bu(\bx;\boldsymbol{H})^\top\boldsymbol{E}\bu(\bx;\boldsymbol{H}). 
  \label{eq:simplicial-Cauchy-Green}
\end{equation}
Thus, on \(\mathrm{GL}^+(N)\), an objective energy density may be written as 
\[ 
  \Psi(\bx,\boldsymbol{H},\boldsymbol{F}) = 
   \widehat{\Psi}(\bx,\boldsymbol{H},\boldsymbol{C}). 
\]

\paragraph{Material isotropy.}
As in the local case, constitutive assumption of material isotropy requires that the stored energy is unaffected by a common rotation of all material directions in the reference configuration.
This requirement may be viewed as invariance under a rotational change of reference description.

To formulate this requirement in the case \(m=N\), consider a nondegenerate reference simplex with edge matrix \(\boldsymbol{H}\) and deformed edge matrix 
\(\boldsymbol{F}\boldsymbol{H}. \)
Let \(\boldsymbol{Q}\in\mathrm{SO}(N)\), and rotate the reference simplex around \(\bx\) according to 
\[ \boldsymbol{H}^* = \boldsymbol{Q}\boldsymbol{H}. \] 
The representation of the deformation gradient relative to the rotated material directions is
\[ \boldsymbol{F}^* = \boldsymbol{F}\boldsymbol{Q}^\top. \] 
Indeed, 
\[ \boldsymbol {F}^*\boldsymbol{H}^* = \boldsymbol{F}\boldsymbol{Q}^\top \boldsymbol{Q}\boldsymbol{H} = \boldsymbol{F}\boldsymbol{H}, \] 
so the pairs 
\[ ( \boldsymbol{H},\boldsymbol{F} ) \qquad\text{and}\qquad ( \boldsymbol{Q}\boldsymbol{H}, \boldsymbol{F}\boldsymbol{Q}^\top ) \] 
describe the same deformed simplex relative to reference simplices that differ only by a rigid rotation.
Material isotropy therefore requires
\[ \Psi( \bx, \boldsymbol{Q}\boldsymbol{H}, \boldsymbol{F}\boldsymbol{Q}^\top) = \Psi\left( \bx, \boldsymbol{H}, \boldsymbol{F} \right),
\qquad  \boldsymbol Q\in\mathrm{SO}(N).\]
In particular, if the energy density does not explicitly depend on \(\boldsymbol{H}\), we recover the condition
\[\Psi(\bx,\boldsymbol{F}\boldsymbol{Q}^\top)=\Psi(\bx,\boldsymbol{F}),\quad \boldsymbol{Q}\in\mathrm{SO}(N).\]

Let us now combine material isotropy with frame indifference.
Under the assumed material rotation, the right Cauchy--Green tensor transforms according to 
\[ \boldsymbol{C}^* = \boldsymbol{Q}\boldsymbol{C}\boldsymbol{Q}^\top. \] 
Consequently, material isotropy requires 
\[ 
\widehat{\Psi}( \bx, \boldsymbol{Q}\boldsymbol{H}, \boldsymbol{Q}\boldsymbol{C}\boldsymbol{Q}^\top) = 
\widehat{\Psi}( \bx, \boldsymbol{H}, \boldsymbol{C}),
\qquad \boldsymbol{Q}\in\mathrm{SO}(N).
\]
Thus, when explicit dependence on the reference simplex is retained, the energy may depend on joint rotational invariants of \(\boldsymbol{C}\) and the geometry of \(\boldsymbol{H}\).
However, for models we primarily have in mind, the dependence on \(\boldsymbol{H}\) enters the energy density only through quantities that are invariant under rotations of the reference configuration, such as 
\[ |\boldsymbol{H}|_{F} \qquad\text{or}\qquad \boldsymbol{H}^\top\boldsymbol{H}. \] 
Indeed, 
\[ |\boldsymbol{Q}\boldsymbol H|_{F} = |\boldsymbol{H}|_{\rm F}, \qquad (\boldsymbol{Q}\boldsymbol{H})^\top (\boldsymbol{Q}\boldsymbol{H}) = \boldsymbol{H}^\top\boldsymbol{H}. \] 
Consequently, an objective and isotropic density may take the form 
\begin{equation} \label{eq:N-dimensional-isotropic-density}
\Psi( \bx, \boldsymbol{H}, \boldsymbol{F} ) = \widehat{\Psi}\bigl( \bx, \boldsymbol{H}^\top\boldsymbol{H}, I_1(\boldsymbol{C}),\ldots,I_N(\boldsymbol{C}) \bigr),\end{equation}
where \(I_1(\boldsymbol{C}),\ldots,I_N(\boldsymbol{C})\) are the principal invariants of the right Cauchy--Green tensor.
Thus, the intrinsic scale and shape of the reference simplex may enter as constitutive parameters without introducing a preferred orientation in the reference configuration.

The invariants are the elementary symmetric functions 
\[ 
  I_k(\boldsymbol{C}) = \sum_{1\leq i_1<\cdots<i_k\leq N} 
  \lambda_{i_1}(\boldsymbol{C})\cdots \lambda_{i_k}(\boldsymbol{C}), \qquad k=1,\ldots,N,
\] 
of the eigenvalues of \(\boldsymbol{C}\).
Equivalently, the deformation dependence of the energy density can be expressed as a symmetric function of the principal stretches
(singular values \(\sigma_i(\boldsymbol{F})\)) of \(\boldsymbol{F}\).
The invariants \(I_k(\boldsymbol{C})\) have a direct geometric interpretation. They aggregate the squared changes of \(k\)-dimensional elements under the linear map \(\boldsymbol{F}\).
In particular, 
\begin{equation} 
    I_1(\boldsymbol{C})=|\boldsymbol{F}|_{F}^2, \qquad 
    I_{N-1}(\boldsymbol{C}) =|\operatorname{cof}\boldsymbol{F}|_{F}^2, \qquad 
    I_N(\boldsymbol{C}) =(\det\boldsymbol{F})^2.
    \label{eq:invariants-minors-general-dimension} 
\end{equation} 
In the mechanically most relevant case \(N=3\), these are precisely the three invariants.
Restricting attention to orientation-preserving deformations \(\boldsymbol{F}\in \mathrm{GL}^+(N)\),
the determinant is determined by the invariants of \(\boldsymbol{C}\). 
In dimension \(N=3\), the isotropic density may therefore equivalently be written as
\begin{equation} 
    \Psi(\bx,\boldsymbol{H}, \boldsymbol{F}) = \widetilde{\Psi}\bigl( \bx, \boldsymbol{H}^\top\boldsymbol{H}, |\boldsymbol F|_{F}, |\operatorname{cof}\boldsymbol{F}|_{F}, \det\boldsymbol{F} \bigr). 
    \label{eq:three-dimensional-isotropic-density}
\end{equation}

\paragraph{Classical and polyconvex constitutive laws.} 
The representation \eqref{eq:N-dimensional-isotropic-density} permits simplex-wise analogues of standard hyperelastic constitutive models. 
The same construction permits polyconvex simplicial densities in the sense of~\cite{Ball1976}. 
In dimension \(N\), such a density has the form 
\begin{equation} 
    \Psi(\bx,\boldsymbol{H},\boldsymbol{F}) = 
    \Psi^{\mathrm{pc}}\bigl( \bx, \boldsymbol{H}^\top\boldsymbol{H}, M_1(\boldsymbol{F}),\ldots,M_N(\boldsymbol{F}) \bigr), 
    \label{eq:polyconvex-simplicial-density} 
\end{equation} 
where \(M_k(\boldsymbol{F})\) denotes the collection of all minors of order \(k\) of 
\(\boldsymbol{F}\), and \(\Psi^{\mathrm{pc}}(\bx,\boldsymbol{H}^\top\boldsymbol{H},\cdot)\) is convex. 
For each fixed \((\bx,\boldsymbol{H}^\top\boldsymbol{H})\), polyconvexity is imposed only with respect to the minors of \(\boldsymbol{F}\). The reference simplex geometry therefore enters parametrically and does not interfere with the convexity requirement.
In three dimensions~\eqref{eq:polyconvex-simplicial-density}  reduces to the familiar representation 
\[ 
  \Psi(\bx,\boldsymbol{H},\boldsymbol{F}) = 
  \Psi^{\mathrm{pc}}\bigl( \bx, \boldsymbol{H}^\top\boldsymbol{H}, \boldsymbol{F}, \operatorname{cof}\boldsymbol{F}, \det\boldsymbol{F} \bigr).
\]
Note that polyconvexity is independent of the preceding constitutive symmetry requirements. 
Accordingly, an objective and isotropic polyconvex density must satisfy the frame-indifference and material-isotropy conditions above in addition to the convexity requirement in the minors of \(\boldsymbol{F}\).

This framework accommodates simplex-wise evaluations of standard constitutive laws from nonlinear solid mechanics, including Saint Venant--Kirchhoff, neo-Hookean, Mooney--Rivlin, and Ogden-type models; see~\cite{Ciarlet1988,Holzapfel2000,Ogden1997}. 
Suitable neo-Hookean and Mooney--Rivlin-type energies provide familiar polyconvex examples, whereas general Saint Venant--Kirchhoff or Ogden-type densities need not be globally polyconvex.

\paragraph{Dependence on the reference simplex $\boldsymbol{H}$.} 
The preceding discussion identifies a natural division of roles between the arguments of the stored-energy density. 
The dependence on \(\boldsymbol{F}\) describes the deformation of the interacting simplex and is subject to the usual constitutive requirements of frame indifference, material symmetry, and, when appropriate, polyconvexity. 
The dependence on \(\boldsymbol{H}^\top\boldsymbol{H}\) describes the intrinsic scale and shape of the reference simplex. 
It may therefore be used to weight interactions according to their size or geometric quality without introducing a preferred orientation in the reference configuration. 
More general dependence on \(\boldsymbol{H}\) may be used to model preferred material directions or other anisotropic effects. 
In that case, the constitutive density should be invariant under the material-symmetry group \(\mathcal{G}\subseteq {\mathrm SO}(N)\), that is,
under the simultaneous transformation
\[ ( \boldsymbol{H},\boldsymbol{F} ) \longmapsto ( \boldsymbol{Q}\boldsymbol{H}, \boldsymbol{F}\boldsymbol{Q}^\top ),
\qquad \text{for every\ }\boldsymbol{Q}\in\mathcal{G}.\] 
Thus, material isotropy does not require the stored-energy density to be independent of the reference simplex, but constrains its dependence on \(\boldsymbol{H}\) through the corresponding material-symmetry condition. 
Polyconvexity, in turn, is imposed on the deformation dependence for each fixed reference-simplex geometry.

\begin{remark}[Relabelling of neighbouring vertices] 
The analytical results below do not require invariance under permutations of the neighbouring vertices. 
If the interacting simplex is regarded as geometrically unlabeled, this invariance may be imposed without affecting the other constitutive assumptions.
Indeed, with the central vertex fixed, let \(\mathfrak{S}_N\) denote the set of \(N\times N\) permutation matrices and define 
\[ 
\Psi_{\rm sym}( \bx,\boldsymbol{H},\boldsymbol{F} ) = \frac{1}{N!} \sum_{\boldsymbol{P}\in\mathfrak{S}_N} 
\Psi( \bx,\boldsymbol{H}\boldsymbol{P},\boldsymbol{F} ). 
\] 
Then \(\Psi_{\rm sym}\) is invariant under relabelling of the neighbouring vertices. 
This symmetrization preserves frame indifference, material isotropy, and polyconvexity in \(\boldsymbol{F}\), as well as any uniform 
growth conditions satisfied by \(\Psi\). 
If the interaction measure is invariant under the right action \(\boldsymbol{H}\mapsto\boldsymbol{H}\boldsymbol{P}\), then a change of variables, together with the invariance of the simplicial deformation gradient under simultaneous relabelling, shows that replacing \(\Psi\) by \(\Psi_{\rm sym}\) does not change the total energy.
\end{remark}

\paragraph{Interaction order and constitutive expressiveness.} 
The number \(m\) of neighbours participating simultaneously in an interaction may be 
viewed as a measure of its constitutive expressiveness.
Indeed, \eqref{eq:Taylor1} implies that a single simplicial gradient 
\(\bD\bu(\bx,\boldsymbol{H})\) with \(m<N\) probes only the restriction of 
\(\boldsymbol{\mathcal{D}}\bu(\bx)\) to an \(m\)-dimensional subspace 
\(\operatorname{range}\boldsymbol{H}\). 
However, a sufficiently rich family of interaction subspaces may collectively determine 
the complete gradient. 
Classical isotropic central-force one-neighbour interactions detect only directional extension, and bond-based models of this type are known to impose restrictions on the resulting local elastic response.
Simultaneous two-neighbour interactions additionally detect coupled directional and 
planar deformation.
At the linearized level, this suggests that two-neighbour interactions may provide enough independent information to distinguish dilatational and distortional response, and hence to represent two independent isotropic elastic moduli.
Finally, when \(m=N\), each nondegenerate simplex determines the complete affine deformation
gradient and an arbitrary hyperelastic density can be applied to it directly.
Turning this kinematic and algebraic hierarchy into a rigorous localization result 
requires further analysis. 
In particular, recovery of the full gradient by averaging does not by itself imply 
recovery of a nonlinear stored-energy density, since averaging and nonlinear constitutive
evaluation do not generally commute.
Determining which local elastic and hyperelastic energies arise as limits for each
interaction order \(m\) is therefore a natural question for the generalized framework 
developed here.

%% file: section_nwspaces.tex
%=====================================================
\section{Radial-\ and shape-weighted nonlocal spaces} 
\label{sec:weighted-NW} 

\subsection{Definition and basic properties}

The definition of the simplicial gradient~\eqref{eq:def_simpl_gradient} treats all 
admissible increment matrices  \(\boldsymbol{H}\) on an equal footing.
However, in applications interactions at different distances need not have the same
strength.
Equally important, \eqref{eq:Taylor1} indicates that finite difference information 
obtained from simplices of different sizes and shapes may have varying reliability, and therefore
does not have to be aggregated with the same weight.
In particular, near-degenerate simplices may amplify finite-difference errors through the 
factor \(\boldsymbol{H}^\dagger\).
We therefore introduce a weighted version of the nonlocal spaces considered 
in~\cite{pedregal_special_2024}, in which the scale and the shape of the interaction are
weighted separately.
To this end, for each \(\boldsymbol{H}\neq {\bf 0}\) we write 
$$
r = |\boldsymbol{H}|_F,\quad 
\boldsymbol{S} = r^{-1} \boldsymbol{H} \in \mathbb{S}^{Nm-1},
$$
where  
$$
\mathbb{S}^{Nm-1} = \{\, \boldsymbol{S} \in \R^{N\times m} \colon |\boldsymbol{S}|_F = 1\,\} 
$$
is the unit sphere in \(\R^{N\times m}\) with respect to Frobenius norm.
The variable \(r\) describes the overall scale of the interaction, whereas 
\(\boldsymbol{S}\) describes the shape and orientation of the associated simplex.
In this notation the simplicial gradient~\eqref{eq:def_simpl_gradient} is defined as follows:
\begin{equation}\label{eq:radial-simplicial-gradient}
    \begin{aligned}
        \boldsymbol{S} &= [\boldsymbol{S}_1,\dots, \boldsymbol{S}_m],\\
        \Delta_{r\boldsymbol{S}} \bu(\bx_0) &= [\bu(\bx_0+r\boldsymbol{S}_1)-\bu(\bx_0),\ldots,\bu(\bx_0+r\boldsymbol{S}_m)-\bu(\bx_0)],\\
        \bD\boldsymbol{u}(\bx_0; r\boldsymbol{S}) &= \frac{1}{r}\Delta_{r\boldsymbol{S}}\bu(\bx_0)\boldsymbol{S}^\dagger.\\
    \end{aligned}
\end{equation}

\paragraph{Radial and shape-orientation measures.} 
The relative importance of different interaction scales and simplex shapes will be modelled using two nonnegative measures. 
Let \(\nu\in\mathcal{M}^+((0,\infty))\) be a nonnegative locally finite Borel measure on
the set of interaction  scales, and let \(\mu\in\mathcal{P}(\mathbb{S}^{Nm-1})\) be a 
Borel probability measure on the set of shapes and orientations.
We refer to \(\nu\) as the \emph{radial measure} and to \(\mu\) as the \emph{shape-orientation measure} (or simply shape measure).
The normalization of \(\mu\) is merely conventional, since its total mass could equally well
be incorporated into \(\nu\). 
No rank condition on \(\boldsymbol{S}\) is imposed at this point. 
Indeed, the Moore--Penrose inverse and the simplicial gradient 
in~\eqref{eq:radial-simplicial-gradient} remain well defined for rank-deficient matrices. 
The possible suppression of nearly degenerate or rank-deficient simplices will instead be encoded through the choice of the shape measure. 
For example, one may impose 
\(
    \mu\left( \left\{ \boldsymbol{S}\colon
    \operatorname{rank}\boldsymbol{S}=m \right\} \right)=1, 
\) 
or a quantitative condition such as 
$$
  \int_{\mathbb{S}^{Nm-1}} |\boldsymbol{S}^\dagger|_{F}^q \,\mathrm{d}\mu(\boldsymbol{S})<\infty 
$$
for a suitable exponent \(q>0\). 
Such assumptions will be introduced only when required.

To account for the boundedness of \(\Omega\),
for  \(\boldsymbol{S} = [\boldsymbol S_1,\ldots,\boldsymbol S_m] \in\mathbb{S}^{Nm-1},\) 
let
\[ 
    r_\O(\bx,\boldsymbol{S}) =  \left\{ r>0: \bx+r\boldsymbol{S}_i\in\O \text{ for every }i=1,\ldots,m \right\}.
\]

\paragraph{Weighted nonlocal  spaces.} 
Let \(1\leq p<\infty\). 
We define the weighted nonlocal seminorm 
\begin{equation} 
    \begin{aligned} 
        [\bu]_{NW_{\nu,\mu}^{1,p}(\O;\R^d)}^p &= 
        \int_\O \int_{\mathbb{S}^{Nm-1}} \int_{r_\O(\bx,\boldsymbol{S})}  | \bD\bu(\bx;r\boldsymbol{S}) |_{F}^p \,\mathrm{d}\nu(r)\,\mathrm{d}\mu(\boldsymbol S)\,\mathrm{d}\bx. 
    \end{aligned} \label{eq:weighted-NW-seminorm} 
\end{equation} 
The corresponding weighted nonlocal  space is 
\begin{equation} 
    \begin{aligned} 
        NW_{\nu,\mu}^{1,p}(\O;\R^d) = 
        \bigl\{ \bu\in L^p(\O;\R^d)\colon 
        [\bu]_{NW_{\nu,\mu}^{1,p}(\O;\R^d)} <\infty \bigr\}, 
    \end{aligned} 
    \label{eq:weighted-NW-space} 
\end{equation} 
equipped with the norm 
\begin{equation} 
    \|\bu\|_{NW_{\nu,\mu}^{1,p}(\O;\R^d)}^p = \|\bu\|_{L^p(\Omega;\mathbb R^d)}^p + 
    [\bu]_{NW_{\nu,\mu}^{1,p}(\O;\R^d)}^p. 
    \label{eq:weighted-NW-norm} 
\end{equation} 

\paragraph{Density-based formulation.} 
In applications, the measures \(\nu\) and \(\mu\) will often be absolutely continuous 
with respect to natural reference measures.
Writing 
\[ 
    \mathrm{d}\nu(r)=\eta(r)\,\frac{\mathrm{d}r}{r}, \qquad 
    \mathrm{d}\mu(\boldsymbol{S}) = \vartheta(\boldsymbol{S})\,\mathrm{d}\sigma(\boldsymbol{S}), 
\] 
where \(\sigma\) denotes surface measure on \(\mathbb{S}^{Nm-1}\), 
the seminorm becomes 
\begin{equation} 
  \begin{aligned} 
    [\bu]_{NW_{\eta,\vartheta}^{1,p}(\O;\R^d)}^p = 
    \int_\Omega \int_{\mathbb{S}^{Nm-1}} \int_{r_\O(\bx,\boldsymbol{S})} 
    |\bD\bu(\bx;r\boldsymbol{S})|_{F}^p \eta(r)\vartheta(\boldsymbol{S}) \,\frac{\mathrm{d}r}{r}\,\mathrm{d}\sigma(\boldsymbol{S})\,\mathrm{d}\bx. 
  \end{aligned} 
  \label{eq:weighted-NW-density-form} 
\end{equation} 
Here \(\eta\) is the radial interaction kernel, while \(\vartheta\) weights the different
simplex shapes and orientations. 
If the weighting is given relative to Lebesgue measure in the increment matrix 
\(\boldsymbol{H}\in\R^{N\times m} \simeq\R^{Nm},\) then the polar-coordinate formula 
\begin{equation} 
 \mathrm{d}\boldsymbol{H} = 
 r^{Nm}\,\frac{\mathrm{d}r}{r}\,\mathrm{d}\sigma(\boldsymbol{S}) 
 \label{eq:matrix-polar-coordinates}
\end{equation}
must be taken into account.
Thus a kernel of the product form 
\[
  \rho(\boldsymbol{H}) = \eta(|\boldsymbol{H}|_{F}) \vartheta\left(\frac{\boldsymbol{H}}{|\boldsymbol{H}|_{F}} \right) 
\]
corresponds to the measure 
\[ 
  \eta(r)\vartheta(\boldsymbol{S}) r^{Nm}\,\frac{\mathrm{d}r}{r}\,\mathrm{d}\sigma(\boldsymbol{S}). 
\]
In the measure formulation~\eqref{eq:weighted-NW-seminorm}, the polar Jacobian and all
additional weights are absorbed into \(\mathrm{d}\nu(r)\,\mathrm{d}\mu(\boldsymbol{S})\).

\paragraph{Relation to the unweighted spaces.} 
The unweighted nonlocal spaces considered in \cite{pedregal_special_2024} are recovered, 
up to an inessential multiplicative constant, by taking the measures induced by Lebesgue
measure in \(\R^{N\times m}\). 
More precisely, on a bounded domain one may take 
\[ 
  \mathrm{d}\nu(r)=r^{Nm}\,\frac{\mathrm{d}r}{r},
  \qquad
  \mathrm{d}\mu(\boldsymbol{S}) = \frac{\mathrm{d}\sigma(\boldsymbol{S})}{\sigma(\mathbb{S}^{Nm-1})}.
\]
Then 
\begin{equation}
  \begin{aligned} 
    \int_{\mathbb{S}^{Nm-1}} \int_{r_\O(\bx,\boldsymbol{S})} 
    |\bD\bu(\bx;r\boldsymbol{S})|_{F}^p r^{Nm} 
    \,\frac{\mathrm{d}r}{r}\,\mathrm{d}\sigma(\boldsymbol{S})
    = \int_{\mathcal{H}(\bx)} |\bD\bu(\bx;\boldsymbol{H})|_{F}^p \,\mathrm{d}\boldsymbol{H}.
  \end{aligned} 
  \label{eq:weighted-unweighted-correspondence}
\end{equation} 
The weighted spaces therefore extend the original construction 
of~\cite{pedregal_special_2024} by allowing interaction scales and simplex shapes to 
contribute nonuniformly.

\paragraph{Examples of radial and shape weights.} 
A finite interaction (peridynamic) horizon \(\delta>0\) can be imposed by requiring 
\[
  \operatorname{spt}\nu\subset(0,\delta]. 
\]
Power-law densities may be used to reproduce fractional-type scalings. 
For \(m=1\), a density behaving as \( \eta(r)\propto r^{p(1-s)} \) reproduces the usual fractional difference-quotient scaling of order \(r^{-sp}\,\mathrm{d}r/r\):
\[
|\bD\bu(\bx;r\boldsymbol{S})|_F^p\eta(r)\,\frac{\mathrm{d}r}{r}
\propto
\frac{|\Delta_{r\boldsymbol{S}}\bu(\bx)\boldsymbol{S}^\dagger|^p}{r^{sp}}\,\frac{\mathrm{d}r}{r}.
\]
For \(m\geq2\), the relation between radial powers and fractional regularity is slightly more involved; see Corollary~\ref{prop:power-law-fractional-embedding} for details.
An isotropic distribution of simplex orientations can be described by a 
rotation-invariant shape measure.
Preferred material directions can instead be incorporated through a nonuniform measure 
\(\mu\).
In particular, discrete shape measures of the form 
\[ \mu = \sum_{\alpha=1}^K a_\alpha\delta_{\boldsymbol{S}_\alpha}, \qquad a_\alpha\geq0, \qquad \sum_{\alpha=1}^K a_\alpha=1, 
\] 
arise naturally in numerical discretizations based on a finite family of simplex 
configurations. 
Near-degenerate shapes may be suppressed for example by choosing a density 
\(\vartheta(\boldsymbol{S}) \propto \det(\boldsymbol{S}^{\top}\boldsymbol{S})^\beta,\) \(\beta> 0\).
Along full-rank shapes, this density tends to zero as the smallest singular value of \(\boldsymbol{S}\) tends to zero and vanishes on rank-deficient matrices.
Conversely, allowing the shape measure to charge rank-deficient matrices provides a 
unified way to include interactions whose effective geometric dimension is lower than \(m\).

\begin{proposition} \label{prop:weighted-NW-Banach} 
  Let \(1\leq p<\infty\), let \(\nu\) be a positive, locally finite Borel measure on 
  \((0,\infty)\), and let \(\mu\) be a Borel probability measure on  \(\mathbb{S}^{Nm-1}\). 
  Then \( NW_{\nu,\mu}^{1,p}(\O;\R^d) \) is a Banach space with respect to the 
  norm~\eqref{eq:weighted-NW-norm}.
\end{proposition} 
\begin{proof}
  Let us put 
   \[
   \chi_\O(\bx,\boldsymbol{S},r) = \begin{cases}
   1, &\qquad \text{if \ } \bx + r \boldsymbol{S}\boldsymbol{e}_i \in \O, \hbox{ for all } i\in \{1,\ldots,m\},\\
   0, &\qquad \text{otherwise}.
   \end{cases}
   \]
   and 
   \[
   \mathrm{d}\Lambda(\bx,\boldsymbol{S},r) = 
   \chi_\O(\bx,\boldsymbol{S},r)\,\mathrm{d}\bx\,\mathrm{d}\mu(\boldsymbol{S})\,\mathrm{d}\nu(r).
   \]

   Let \(\{\bu_j\}\) be Cauchy in \(NW_{\nu,\mu}^{1,p}(\O;\R^d)\). 
   Then, for some \(\bu\in L^p(\O;\R^d)\) and \(\boldsymbol{G}\in L^p(\Lambda;\R^{d\times N})\), 
   \[\begin{aligned} 
    \bu_j&\to\bu, &\quad&\text{in }L^p(\O;\R^d), \qquad\text{and}\\
    \bD\bu_j&\to\boldsymbol{G}, &\quad&\text{in }L^p(\Lambda;\R^{d\times N}).
   \end{aligned}\]
   Passing to subsequences, both convergences hold almost everywhere in their 
   respective measure spaces. 
   For every fixed admissible \((\boldsymbol{S},r)\), translations preserve 
   Lebesgue-null sets. 
   Fubini's theorem therefore gives, for \(\Lambda\)-almost every 
   \((\bx,\boldsymbol{S},r)\), 
   \[ 
     \Delta_{r\boldsymbol{S}}\bu_j(\bx) \longrightarrow \Delta_{r\boldsymbol{S}}\bu(\bx).
   \] 
   Since \(\boldsymbol{S}^\dagger\) is fixed at this point, 
   \[ 
     \bD\bu_j(\bx;r\boldsymbol{S}) \longrightarrow \bD\bu(\bx;r\boldsymbol{S}). 
   \]
   Therefore \[ \boldsymbol{G}=\bD\bu \quad\text{\(\Lambda\)-almost everywhere}. \] 
   It follows that \(\bu\in NW_{\nu,\mu}^{1,p}\) and that \(\bu_j\to\bu\) 
   in the norm~\eqref{eq:weighted-NW-norm}.
\end{proof}

\begin{corollary} 
  Let \(1<p<\infty\). Then \(NW_{\nu,\mu}^{1,p}(\O;\R^d)\) is reflexive.
  For \(p=2\), it is a Hilbert space.
\end{corollary}
\begin{proof} 
  The map \(\bu \longmapsto (\bu,\bD\bu) \) identifies 
  \(NW_{\nu,\mu}^{1,p}(\O;\R^d)\) with a closed subspace of 
  \( L^p(\O;\R^d) \times L^p(\Lambda;\R^{d\times N}). \)
  The conclusion follows from the reflexivity of \(L^p\) for \(1<p<\infty\). 
  For \(p=2\), the product space is Hilbert. 
\end{proof}

%=====================================================
\subsection{Edge marginals and comparison with two-point energies}
\label{subsec:edge-marginals}

A single simplicial gradient contains finite-difference information
along each edge of the underlying simplex. We now make this
observation quantitative by introducing the corresponding
one-edge marginals of the radial and shape-orientation measures.

For the results in this subsection, we assume that
\begin{equation}
    \mu\left(
        \left\{
            \boldsymbol{S}\in\mathbb{S}^{Nm-1}:
            \operatorname{rank}\boldsymbol{S}=m
        \right\}
    \right)
    =
    1.
    \label{eq:full-column-rank-shapes}
\end{equation}
In particular,
\(
    \boldsymbol{S}^\dagger\boldsymbol{S}
    =
    \boldsymbol{I}_m,
\)
for \(\mu\)-almost every \(\boldsymbol{S}\).

For \(\bx\in\O\), \(i\in\{1,\ldots,m\}\), and every Borel set
\(B\subset\mathbb{R}^N\), let
\begin{equation}
    \kappa_{\bx,i}(B)
    =
    \int_{\mathbb{S}^{Nm-1}}
    \int_{r_\Omega(\bx,\boldsymbol S)}
    \chi_B(
        r\boldsymbol S\boldsymbol e_i
    )
    \,\mathrm d\nu(r)\,\mathrm d\mu(\boldsymbol S).
    \label{eq:marginal-measure}
\end{equation}
In other words, \(\kappa_{\bx,i}\) is the push-forward of the
restriction of \(\nu\otimes\mu\) to the admissible simplex
configurations under the map
\(
    (r,\boldsymbol S)
    \longmapsto
    r\boldsymbol S\boldsymbol e_i.
\)
We also set
\begin{equation}
    \kappa_{\bx}
    =
    \sum_{i=1}^m\kappa_{\bx,i}.
    \label{eq:total-marginal-measure}
\end{equation}
We refer to \(\kappa_{\bx,i}\) as the \(i\)-th edge marginal, and to
\(\kappa_{\bx}\) as the total edge marginal (of $\nu\otimes\mu$).

The edge marginals give rise to the two-point energy
\begin{equation}
\begin{aligned}
    \mathcal{Q}_{\nu,\mu}(\bu)
    &=
    \int_\O
    \int_{\R^N}
    \frac{
        |\bu(\bx+\boldsymbol h)-\bu(\bx)|^p
    }{
        |\boldsymbol h|^p
    }
    \,\mathrm{d}\kappa_{\bx}(\boldsymbol{h})\,\mathrm{d}\bx.
\end{aligned}
\label{eq:edge-marginal-energy}
\end{equation}
Since \(\kappa_{\bx}\) is supported in \(\O-\bx\), the value
\(\bu(\bx+\boldsymbol{h})\) appearing in
\eqref{eq:edge-marginal-energy} is always well defined up to the
usual choice of an almost-everywhere representative.

\begin{proposition}[Comparison with the edge-marginal energy]
\label{prop:edge-marginal-comparison}
Let \(1\leq p<\infty\), and suppose that
\eqref{eq:full-column-rank-shapes} holds. Then, for every
\(\bu\in NW_{\nu,\mu}^{1,p}(\O;\R^d)\),
\begin{equation}
    \mathcal{Q}_{\nu,\mu}(\bu)
    \leq
    m
    [\bu]_{NW_{\nu,\mu}^{1,p}(\O;\R^d)}^p.
    \label{eq:edge-marginal-comparison}
\end{equation}
\end{proposition}

\begin{proof}
Let
\(
    \boldsymbol{S}
    =
    [
        \boldsymbol{S}_1,\ldots,\boldsymbol{S}_m
    ]
\)
have full column rank. For every admissible
\((\bx,\boldsymbol{S},r)\) and every \(i\in\{1,\ldots,m\}\), we have
\[
    \bD\bu(\bx;r\boldsymbol{S})
    (r\boldsymbol{S}\boldsymbol{e}_i)
    =
    \frac{1}{r}
    \Delta_{r\boldsymbol{S}}\bu(\bx)
    \boldsymbol{S}^\dagger
    (r\boldsymbol{S}\boldsymbol{e}_i)
    =
    \Delta_{r\boldsymbol{S}}\bu(\bx)
    \boldsymbol{e}_i
    =
    \bu(\bx+r\boldsymbol{S}_i)-\bu(\bx),
    \]
and therefore
\[
    \frac{
        |\bu(\bx+r\boldsymbol{S}_i)-\bu(\bx)|
    }{
        |r\boldsymbol{S}_i|
    }
    \leq
    |\bD\bu(\bx;r\boldsymbol{S})|_{F}.
\]
Raising this inequality to the power \(p\), summing over
\(i=1,\ldots,m\), and integrating gives
\[\begin{aligned}
    \mathcal{Q}_{\nu,\mu}(\bu)=&\int_\O
    \int_{\mathbb{S}^{Nm-1}}
    \int_{r_\O(\bx,\boldsymbol{S})}
    \sum_{i=1}^m
    \frac{
        |\bu(\bx+r\boldsymbol{S}_i)-\bu(\bx)|^p
    }{
        |r\boldsymbol{S}_i|^p
    }
    \,\mathrm{d}\nu(r)\,\mathrm{d}\mu(\boldsymbol{S})\,\mathrm{d}\bx
    \\
    &\leq
    m
    \int_\O
    \int_{\mathbb{S}^{Nm-1}}
    \int_{r_\O(\bx,\boldsymbol{S})}
    |\bD\bu(\bx;r\boldsymbol{S})|_{F}^p
    \,\mathrm{d}\nu(r)\,\mathrm{d}\mu(\boldsymbol{S})\,\mathrm{d}\bx
    =m
    [\bu]_{NW_{\nu,\mu}^{1,p}(\O;\R^d)}^p. 
    \qedhere
\end{aligned}\]
\end{proof}
The total edge marginal may be interpreted as the effective distribution of constituent bonds generated by the multipoint interaction family. 
Proposition~\ref{prop:edge-marginal-comparison} shows that every simplicial interaction controls the corresponding collection of two-point difference quotients.
To obtain explicit lower bounds on the edge marginals for isotropic shape measures, which are subsequently utilized to establish non-local Poincaré inequalities and continuous injections in standard fractional Sobolev spaces, 
we shall use the following uniform interior-thickness condition on \(\O\):
\begin{quote}
there exist \(c_\O,\delta_\O>0\) such that, for every \(\bx\in\O\) and every \(0<t<\delta_\O\), 
\begin{equation} 
    \left| (\O-\bx) \cap \bigl(B_{2t}(\boldsymbol{0})\setminus B_t(\boldsymbol{0})\bigr) \right| \geq 
    c_\O t^N.
    \label{eq:uniform-interior-thickness} 
\end{equation}
\end{quote}
Condition~\eqref{eq:uniform-interior-thickness} follows from the uniform interior cone property and is therefore satisfied by every bounded Lipschitz domain; see for example \cite{Grisvard1985}.

\begin{lemma}[Edge marginals for isotropic power-law weights] 
  \label{lem:power-law-edge-marginals} 
  Suppose $\O$ is a bounded domain. 
  Suppose that \(\mu\) is normalized surface measure, and 
  \[ 
    \mathrm d\nu(r) = \eta(r)\frac{\mathrm dr}{r}, \qquad 
    \eta(r)\geq c_\eta r^a \quad\text{for }0<r<\delta, a>0. 
  \] 
  \begin{enumerate}
    \item For \(m=1\),
    \begin{equation}\label{eq:scale-marginal-lower-bound-m1}
      \mathrm{d}\kappa_{\bx}(\boldsymbol{h}) \geq 
      c \chi_{ (\O-\bx)\cap B_\delta(\boldsymbol{0}) }(\boldsymbol{h}) 
      |\boldsymbol{h}|^{a-N} \,\mathrm{d}\boldsymbol{h}.
    \end{equation}
    \item Let \(2\leq m\leq N\), and suppose that~\eqref{eq:uniform-interior-thickness} holds. 
    Set 
    \begin{equation}\label{eq:t-star-lemma}
    t_* = \min\left\{ \frac{\delta}{2\sqrt m}, \delta_\O \right\}.
    \end{equation}
    Then there exists \(c_*>0\), depending only on \(a,N,m,c_\eta,c_\O\), such that, 
    for every \(t\in(0,t_*)\) and almost every \(\bx\in\O\), 
    \begin{equation} 
      \mathrm{d}\kappa_{\bx,1}(\boldsymbol{h}) \geq c_* t^{a-N} 
      \chi_{ (\O-\bx)\cap B_t(\boldsymbol{0}) }(\boldsymbol{h}) 
      \,\mathrm{d}\boldsymbol{h}.
      \label{eq:scale-dependent-marginal-lower-bound}
    \end{equation} 
    Consequently, the same lower bound holds with \(\kappa_{\bx}\) in place of 
    \(\kappa_{\bx,1}\).
  \end{enumerate} 
\end{lemma}
\begin{proof}
  {\bfseries Case \(m=1\).}
  For every Borel set \(B\subset\R^N\), the polar-coordinate formula gives 
  \[ 
  \begin{aligned} \kappa_{\bx,1}(B) &= 
    \frac{1}{\sigma(\mathbb S^{N-1})} \int_{B\cap(\Omega-\bx)} \eta(|\boldsymbol h|) |\boldsymbol h|^{-N} \,\mathrm d\boldsymbol h. 
  \end{aligned} \] 
  Consequently, 
  \[ \mathrm d\kappa_{\bx,1}(\boldsymbol h) \geq \frac{c_\eta}{\sigma(\mathbb S^{N-1})} \chi_{ (\Omega-\bx)\cap B_\delta(\boldsymbol 0) }(\boldsymbol h) |\boldsymbol h|^{a-N} \,\mathrm d\boldsymbol h. 
  \]

  \medskip
  {\bfseries Case \(2\leq m\leq N\).}
  We write \(\boldsymbol{H} = [\boldsymbol{h},\boldsymbol{Z}]\), 
  where \(\boldsymbol{h}\in\R^N\) and \(\boldsymbol{Z}\in\R^{N\times(m-1)}\).
  The measure induced on increment matrices has density 
  \(\sigma(\mathbb{S}^{Nm-1})^{-1} \eta(|\boldsymbol{H}|_{F}) |\boldsymbol{H}|_{F}^{-Nm}\). Consequently, for every Borel set 
  \(B\subset\R^N\), 
  \begin{equation}
    %\begin{aligned} 
      \kappa_{\bx,1}(B)\geq\frac{c_\eta}{\sigma(\mathbb{S}^{Nm-1})} 
      \int_{B\cap(\O-\bx)\cap B_\delta(\boldsymbol{0})} 
      \int_{\mathcal{Z}_\delta(\bx,\boldsymbol{h})} \bigl( |\boldsymbol{h}|^2+|\boldsymbol{Z}|_{\rm F}^2 \bigr)^{(a-Nm)/2} \,\mathrm{d}\boldsymbol{Z}\,\mathrm{d}\boldsymbol{h}, \label{eq:marginal-power-law-preestimate} 
    %\end{aligned}
  \end{equation} 
  where 
  \[ 
    \mathcal{Z}_\delta(\bx,\boldsymbol{h}) =  \left\{ \boldsymbol{Z}\in\R^{N\times(m-1)}: 
    [\boldsymbol{h},\boldsymbol{Z}]\in\mathcal{H}(\bx), \quad 
    |\boldsymbol{h}|^2+|\boldsymbol{Z}|_{F}^2<\delta^2 \right\}. 
  \] 
  Fix \(t\in(0,t_*)\), and define 
  \[ 
    A_t(\bx) = (\O-\bx) \cap 
    \left( B_{2t}(\boldsymbol{0}) \setminus B_t(\boldsymbol{0}) \right). 
  \] 
  Due to~\eqref{eq:uniform-interior-thickness}, 
  \( |A_t(\bx)| \geq c_\O t^N. \) 
  If \(|\boldsymbol{h}|\leq t\), then 
  \( A_t(\bx)^{m-1} \subset \mathcal{Z}_\delta(\bx,\boldsymbol{h}). \) 
  Moreover, for every  \(\boldsymbol{Z}\in A_t(\bx)^{m-1}\), 
  we have the inequalities
  \[ (m-1)t^2 \leq |\boldsymbol{h}|^2+|\boldsymbol{Z}|_{F}^2 < 4mt^2. \] 
  Hence, 
  \[ \bigl( |\boldsymbol{h}|^2+|\boldsymbol{Z}|_{F}^2 \bigr)^{(a-Nm)/2} 
     \geq c_{a,N,m}t^{a-Nm},
  \] 
  where 
  \[ c_{a,N,m} = \begin{cases} (4m)^{(a-Nm)/2}, & a<Nm,\\[1mm] 
                               (m-1)^{(a-Nm)/2}, & a\geq Nm. \end{cases} 
  \]
  Therefore
  \[\int_{\mathcal{Z}_\delta(\bx,\boldsymbol{h})} 
  \bigl( |\boldsymbol{h}|^2+|\boldsymbol{Z}|_{F}^2 \bigr)^{(a-Nm)/2} 
  \,\mathrm{d}\boldsymbol{Z}
  \geq c_{a,N,m}t^{a-Nm} |A_t(\bx)|^{m-1}
  \geq c_\Omega^{m-1}c_{a,N,m}t^{a-N}.
  \]
  Inserting this estimate into~\eqref{eq:marginal-power-law-preestimate} yields~\eqref{eq:scale-dependent-marginal-lower-bound}, 
  with 
  \[ c_* = \frac{ c_\eta c_\O^{m-1}c_{a,N,m} }{ \sigma(\mathbb{S}^{Nm-1}) }. \] 
  Since \(\kappa_{\bx}\geq\kappa_{\bx,1}\), the final assertion follows.
\end{proof}

%===================================================== 
\subsection{Nonlocal Poincar{\'e} inequalities} \label{subsec:nonlocal-Poincare} 
The seminorm \( [\cdot]_{NW_{\nu,\mu}^{1,p}} \) vanishes on constant functions. 
To obtain control of the complete \(L^p\)-norm, one must therefore impose a normalization 
or a Dirichlet-type condition. 
Let \(\O_D\subset\O\) be measurable with \(|\O_D|>0, \) and define 
\begin{equation} 
  \begin{aligned} 
    NW_{\nu,\mu,D}^{1,p}(\O;\R^d) = \Bigl\{ \bu\in NW_{\nu,\mu}^{1,p}(\O;\R^d)
    \colon \bu={\bf 0} \text{ almost everywhere on }\O_D \Bigr\}. 
  \end{aligned} 
  \label{eq:homogeneous-nonlocal-Dirichlet-space}
\end{equation} 
In peridynamic applications, \(\O_D\) is often chosen to be a boundary layer whose
thickness is related to the nonlocal interaction horizon.
We recall that finite-range two-point energies satisfy a nonlocal Poincar{\'e} inequality.
More precisely, if \(\O\) is bounded and connected and \(\delta>0\), then 
\begin{equation}
  \begin{aligned}
  \|\bu-(\bu)_\O\|_{L^p(\O;\R^d)}^p &\leq C \int_\O\int_{\O\cap B_\delta(\bx)} 
  |\bu(\bx')-\bu(\bx)|^p \,\mathrm{d}\bx'\,\mathrm{d}\bx, \qquad\text{where}\\
  (\bu)_\O &= \frac{1}{|\O|} \int_\O\bu(\bx)\,\mathrm{d}\bx.
  \end{aligned}
  \label{eq:finite-range-Poincare-Wirtinger} 
\end{equation}
This is a standard nonlocal Poincar{\'e} inequality; see, for example, \cite{BCM14, Ponce2004} or~\cite{BellidoMoraCorralSchonberger}.

If \(\bu={\bf 0}\) almost everywhere on \(\O_D\), then 
\[
  |\O_D|\,|(\bu)_\O|^p = \int_{\O_D} |(\bu)_\O-\bu(\bx)|^p \,\mathrm{d}\bx 
  \leq \|\bu-(\bu)_\O\|_{L^p(\O;\R^d)}^p.
\] 
Consequently, 
\begin{equation}\label{eq:Poincare-from-OD}
  \|\bu\|_{L^p(\O;\R^d)}^p
  \leq 2^{p-1} \left( \|\bu-(\bu)_\O\|_{L^p(\O;\R^d)}^p + |\O|\,|(\bu)_\O|^p \right)
  \leq 2^{p-1} \left( 1+\frac{|\O|}{|\O_D|} \right) \|\bu-(\bu)_\O\|_{L^p(\O;\R^d)}^p. 
\end{equation} 
Thus \eqref{eq:finite-range-Poincare-Wirtinger} also gives a Poincar{\'e} inequality on 
the homogeneous Dirichlet class.

\begin{proposition}[Nonlocal Poincar{\'e} inequality] 
  \label{prop:NW-Poincare} 
  Let \(1\leq p<\infty\), let \(\O\subset\R^N\) be a bounded connected open set, and 
  suppose that~\eqref{eq:full-column-rank-shapes} holds. 
  Assume that there exist \(c_0,\delta_0>0\) such that, for almost every \(\bx\in\O\), 
  \begin{equation} 
    \mathrm{d}\kappa_{\bx}(\boldsymbol{h}) \geq 
    c_0 \chi_{ (\O-\bx)\cap B_{\delta_0}(\boldsymbol{0}) }(\boldsymbol{h}) 
    |\boldsymbol{h}|^p \,\mathrm{d}\boldsymbol{h}. 
    \label{eq:Poincare-marginal-lower-bound} 
  \end{equation} 
  Let \(\O_D\subset\O\) be measurable with \(|\O_D|>0\). 
  Then there exists \( C= C( p,\O,\O_D,c_0,\delta_0,m )>0 \)
  such that 
  \begin{equation} 
    \|\bu\|_{L^p(\O;\R^d)}^p \leq C [\bu]_{NW_{\nu,\mu}^{1,p}(\O;\R^d)}^p,
    \qquad \forall\bu\in NW_{\nu,\mu,D}^{1,p}(\O;\R^d).
    \label{eq:NW-Poincare} 
  \end{equation} 
\end{proposition}
\begin{proof} 
  Thanks to Proposition~\ref{prop:edge-marginal-comparison} and the bound~\eqref{eq:Poincare-marginal-lower-bound}, we have the inequalities 
  \[ 
    c_0 \int_\O \int_{\O\cap B_{\delta_0}(\bx)} |\bu(\bx')-\bu(\bx)|^p 
    \,\mathrm{d}\bx'\,\mathrm{d}\bx
    \leq \mathcal Q_{\nu,\mu}(\bu) \leq m [\bu]_{NW_{\nu,\mu}^{1,p}}^p. 
  \] 
  The finite-range nonlocal Poincar{\'e} inequality~\eqref{eq:finite-range-Poincare-Wirtinger}, together with the fact that \(\bu={\bf 0}\) almost everywhere on \(\Omega_D\) and~\eqref{eq:Poincare-from-OD}, yields~\eqref{eq:NW-Poincare}. 
\end{proof}

The marginal lower bound thus acts as a connectivity condition on the interaction family: together with the nonlocal Dirichlet constraint, it rules out nonconstant zero-energy modes and yields coercivity of the model.
Combining Proposition~\ref{prop:NW-Poincare} with the explicit edge-marginal estimates from Lemma~\ref{lem:power-law-edge-marginals} gives the following corollaries.

\begin{corollary}[Radial criterion for \(m=1\)] 
  \label{cor:Poincare-power-law-m1} 
  Let \(m=1\), let \(\mu\) be normalized surface measure on \(\mathbb{S}^{N-1}\), 
  and let 
  \[ \mathrm{d}\nu(r) = \eta(r)\frac{\mathrm{d}r}{r}. \] 
  Suppose that, for some \(0<a\leq N+p\), \(c_\eta>0\), \(\delta>0\), we have 
  the estimate \(\eta(r)\geq c_\eta r^a\), for all \( r \in (0,\delta).\)
  Then \eqref{eq:Poincare-marginal-lower-bound} holds, and so does the conclusion of Proposition~\ref{prop:NW-Poincare}.
\end{corollary}
\begin{proof}
  We apply Lemma~\ref{lem:power-law-edge-marginals} to get the 
  lower bound~\eqref{eq:scale-marginal-lower-bound-m1}.
  Since \(a\leq N+p\), we have, for \(0<|\boldsymbol h|<\delta\), 
  \[ 
  |\boldsymbol h|^{a-N} = |\boldsymbol h|^p |\boldsymbol h|^{a-N-p} \geq 
  \min\{1,\delta^{a-N-p}\} |\boldsymbol h|^p. 
  \] 
\end{proof}

\begin{corollary}[Radial criterion for \(m\geq 2\)] 
  \label{cor:Poincare-power-law-mge2} 
  Let \(2\leq m\leq N\), let \(\O\) satisfy~\eqref{eq:uniform-interior-thickness}, and let 
  \(\mu\) be the normalized surface measure on \(\mathbb{S}^{Nm-1}\).
  Suppose that 
  \[ \mathrm {d}\nu(r) = \eta(r)\frac{\mathrm dr}{r} \] and that, for some \(a>0\) and 
  \(c_\eta>0\), \(\delta>0\), we have the estimate 
  \( \eta(r)\geq c_\eta r^a\), for all \(r \in (0,\delta). \) 
   Then \eqref{eq:Poincare-marginal-lower-bound} holds, and so does the conclusion of Proposition~\ref{prop:NW-Poincare}.
\end{corollary}
\begin{proof} 
  Let us fix an arbitrary \(t\in(0,t_*)\), where \(t_*>0\) is given by~\eqref{eq:t-star-lemma}.
  We apply Lemma~\ref{lem:power-law-edge-marginals} to get the estimate~\eqref{eq:scale-dependent-marginal-lower-bound}.
  Note that \(1 \geq t^{-p}|\boldsymbol{h}|^p\),  for \(|\boldsymbol h|<t\).
  For such \(\boldsymbol{h}\),
  we can rewrite~\eqref{eq:scale-dependent-marginal-lower-bound} as
  \[ 
    \mathrm{d}\kappa_{\bx}(\boldsymbol{h}) \geq c_*t^{a-N-p} 
    \chi_{ (\O-\bx)\cap B_{t}(\boldsymbol{0}) }(\boldsymbol{h}) |\boldsymbol{h}|^p 
    \,\mathrm {d}\boldsymbol{h}.
  \qedhere\]
\end{proof}

\subsection{Embedding into fractional Sobolev spaces.} 

The edge-marginal comparison also yields fractional regularity when the total edge marginal dominates a suitable fractional kernel.

\begin{theorem}[Fractional Sobolev embedding]
  \label{thm:fractional-Sobolev-embedding} 
  Let \(1\leq p<\infty\) and \(s\in(0,1)\), and suppose that~\eqref{eq:full-column-rank-shapes} holds. 
  Assume that there exist \(c_0,\delta_0>0\) such that, for almost every \(\bx\in\Omega\), 
  \begin{equation} 
    \mathrm{d}\kappa_{\bx}(\boldsymbol{h}) \geq 
    c_0 \chi_{(\O-\bx)\cap B_{\delta_0}(\boldsymbol{0})}(\boldsymbol{h}) 
    |\boldsymbol{h}|^{p-N-sp}\,\mathrm{d}\boldsymbol{h}. 
    \label{eq:local-marginal-lower-bound} 
  \end{equation}
  Then 
  \begin{equation} 
    NW_{\nu,\mu}^{1,p}(\O;\R^d) \hookrightarrow W^{s,p}(\O;\R^d)
    \label{eq:fractional-Sobolev-embedding} 
  \end{equation}
  continuously. 
  More precisely, there exists a constant \(C=C(p,s,\Omega,\delta_0,c_0,m)>0\) such that
  \begin{equation} 
    \|\bu\|_{W^{s,p}(\O;\R^d)}^p \leq 
    C \left( \|\bu\|_{L^p(\O;\R^d)}^p + [\bu]_{NW_{\nu,\mu}^{1,p}(\O;\R^d)}^p \right).
    \label{eq:fractional-Sobolev-estimate}
  \end{equation} 
\end{theorem}
\begin{proof}
  Utilizing~\eqref{eq:local-marginal-lower-bound}  and Proposition~\ref{prop:edge-marginal-comparison}, we get the estimates
  \[\begin{aligned}
    c_0 \int_\O \int_{(\O-\bx)\cap B_{\delta_0}(\boldsymbol{0})} \frac{ |\bu(\bx+\boldsymbol{h})-\bu(\bx)|^p }{ |\boldsymbol{h}|^{N+sp} } \,\mathrm{d}\boldsymbol{h}\,\mathrm{d}\bx 
    &\leq     
    \int_\O\int_{\R^N} \frac{ |\bu(\bx+\boldsymbol{h})-\bu(\bx)|^p }{ |\boldsymbol{h}|^p } 
    \,\mathrm{d}\kappa_{\bx}(\boldsymbol{h})\,\mathrm{d}\bx 
    \\&\leq m [\bu]_{NW_{\nu,\mu}^{1,p}(\O;\R^d)}^p.
  \end{aligned}\]
 To conclude the proof, we need to estimate the contributions to Gagliardo seminorm
 corresponding to \(|\boldsymbol{h}|\geq \delta_0\).
 To this end, using the inequality \(|a-b|^p \leq 2^{p-1}\bigl(|a|^p+|b|^p\bigr),\)
we obtain 
\[\begin{aligned} 
    \int_\O \int_{(\O-\bx)\cap(\R^N\setminus B_{\delta_0}(\boldsymbol{0}))} 
    \frac{ |\bu(\bx+\boldsymbol{h})-\bu(\bx)|^p }{ |\boldsymbol{h}|^{N+sp} } 
    \,\mathrm{d}\boldsymbol{h}\,\mathrm{d}\bx 
    &\leq 2^{p-1}\delta_0^{-N-sp} \int_\O\int_{\O-\bx} \left( |\bu(\bx+\boldsymbol{h})|^p+|\bu(\bx)|^p \right) \,\mathrm{d}\boldsymbol{h}\,\mathrm{d}\bx 
    \\ &\leq 2^p|\O|\delta_0^{-N-sp} \|\bu\|_{L^p(\O;\R^d)}^p. 
\end{aligned}\]
 Combining the two estimates gives 
 \[ 
   [\bu]_{W^{s,p}(\O;\R^d)}^p \leq 
   \frac{m}{c_0} [\bu]_{NW_{\nu,\mu}^{1,p}}^p + 
   2^p|\O|\delta_0^{-N-sp} \|\bu\|_{L^p(\O;\R^d)}^p. 
 \] 
 Adding \(\|\bu\|_{L^p(\O;\R^d)}^p\) to this inequality concludes the proof.
\end{proof}

The stronger small-scale marginal bound has a distinct analytical role. 
It ensures that the interaction family samples sufficiently many short edges to yield compactness, which in turn allows the existence theory to accommodate nonconvex energy densities.
The following result is a standard consequence of the characterization of fractional Sobolev spaces by difference quotients; see, for example, \cite{MR2944369}.
Related anisotropic fractional seminorms are considered in~\cite{Ludwig2014}.

\begin{corollary}[Power-law radial weights and fractional regularity, \(m=1\)] 
    \label{prop:power-law-bond-embedding} 
    Let \(\O\subset\R^N\) be a bounded open set, let \(1\leq p<\infty\), and let \(m=1\).
    Suppose that \(\mu\) is the normalized surface measure on \(\mathbb{S}^{N-1}\), and that \(\mathrm{d}\nu(r)=\eta(r)\frac{\mathrm{d}r}{r}.\)
    Assume that, for some \(a\in(0,p)\), \(\delta>0\), and \(c_\eta>0\), we have the inequality 
    \[ 
      \eta(r)\geq c_\eta r^a, \qquad 0<r<\delta.
    \] 
    Then \[NW_{\nu,\mu}^{1,p}(\O;\R^d) \hookrightarrow W^{1-a/p,p}(\O;\R^d) \] continuously.
\end{corollary}
\begin{proof}
  We apply Lemma~\ref{lem:power-law-edge-marginals} to get the 
  lower bound~\eqref{eq:scale-marginal-lower-bound-m1}.
  We then put \(s = 1-\frac{a}{p} \in (0,1)\) so that \(a-N=p-N-sp\).
  The conclusion follows from Theorem~\ref{thm:fractional-Sobolev-embedding}.
\end{proof}

\begin{corollary}[Power-law radial weights and fractional regularity, \(m\geq 2\)] 
    \label{prop:power-law-fractional-embedding} 
    Let \(\O\subset\R^N\) be a bounded domain satisfying~\eqref{eq:uniform-interior-thickness}, let \(1\leq p<\infty\), 
    and let \(2\leq m\leq N\).
    Suppose that \(\mu\) is the normalized surface measure on \(\mathbb{S}^{Nm-1}\),
    and \(d\nu(r)=\eta(r)\frac{dr}{r}.\) 
    Assume that, for some \(a>0\), \(\delta>0\), and \(c_\eta>0\), 
    \begin{equation} 
      \eta(r)\geq c_\eta r^a, \qquad 0<r<\delta. 
      \label{eq:radial-power-lower-bound} 
    \end{equation} 
    Then the following statements hold: 
    \begin{enumerate} 
        \item If \(0<a<\min\{N,p\}\), then 
        \begin{equation} 
            NW_{\nu,\mu}^{1,p}(\O;\R^d) \hookrightarrow W^{1-a/p,p}(\O;\R^d)
            \label{eq:subcritical-power-embedding}
        \end{equation} 
        continuously. 
        \item If \(a\geq N\) and \(p>N\), then 
        \begin{equation} 
            NW_{\nu,\mu}^{1,p}(\O;\R^d) \hookrightarrow W^{1-N/p,p}(\O;\R^d)
            \label{eq:saturated-power-embedding}
        \end{equation} 
        continuously.
    \end{enumerate}
\end{corollary}
\begin{proof}
  \emph{Suppose first that \(0<a<\min\{N,p\}. \)}
  For every \(0<|\boldsymbol{h}|< t_*/2\), where \(t_*>0\) is given by~\eqref{eq:t-star-lemma},
  we use Lemma~\ref{lem:power-law-edge-marginals} to get the estimate~\eqref{eq:scale-dependent-marginal-lower-bound}.
  We use this estimate with \(t=2|\boldsymbol{h}|\), so that \(\chi_{ B_{t}(\boldsymbol{0}) }(\boldsymbol{h})\equiv 1\) to obtain the lower bound
  \[ 
    \mathrm{d}\kappa_{\bx}(\boldsymbol{h}) \geq 
    c_*2^{a-N} \chi_{ (\O-\bx)\cap B_{t_*/2}(\boldsymbol{0}) }(\boldsymbol{h}) |\boldsymbol{h}|^{a-N} \,\mathrm{d}\boldsymbol{h}. 
  \] 
  Setting \(s=1-\frac{a}{p}\in(0,1),\) we have \(a-N=p-N-sp,\)
  and the first assertion follows from Theorem~\ref{thm:fractional-Sobolev-embedding}.

  \medskip
  \emph{Suppose next that \(a\geq N\) and \(p>N\).}
  Fix an arbitrary \(t\in(0,t_*)\), where as before \(t_*>0\) is given by~\eqref{eq:t-star-lemma}.
  We again appeal to Lemma~\ref{lem:power-law-edge-marginals} to get the estimate~\eqref{eq:scale-dependent-marginal-lower-bound}.
  Setting \(s=1-\frac{N}{p}\in(0,1),\) we have \(p-N-sp=0\), and the second assertion follows from Theorem~\ref{thm:fractional-Sobolev-embedding}.
\end{proof}

The second case in Corollary~\ref{prop:power-law-fractional-embedding} exhibits a saturation effect: once \(a\geq N\), 
integration over the unused simplex edges produces a one-edge marginal that locally dominates the Lebesgue measure.
The fractional order guaranteed by this argument therefore remains \(1-N/p\), independently of any further increase in \(a\).
The reason for this is that for \(m\geq2\), a short selected edge may belong to many simplices, whose remaining edges stay at a fixed nonzero scale. Integration over these unused edges can therefore produce a one-edge marginal that locally dominates Lebesgue measure, even when the radial weight decays more rapidly. 

The following simple corollary generalizes the observation made in~\cite{pedregal_special_2024} to the case \(1\leq m \leq N\).
Since constant normalizations of the interaction measure do not affect the resulting spaces or embeddings, we suppress the factor \(\sigma(\mathbb{S}^{Nm-1})^{-1}\) in the unweighted formulas below.

\begin{corollary}[Unweighted spaces] 
    \label{cor:Pedregal-from-power-law} 
    Let \(\O\subset\R^N\) be a bounded open set.
    For unweighted integration in \(\R^{N\times m}\), see~\eqref{eq:weighted-unweighted-correspondence},
    whenever \(p>N\), \(NW^{1,p}(\O;\R^d) \hookrightarrow W^{1-N/p,p}(\O;\R^d)\) continuously. 
    For \(m=N\), this recovers the embedding established in~\cite{pedregal_special_2024}.
\end{corollary}
\begin{proof}
  We note that in this case we have \(\eta(r)=r^{Nm}.\)
  If \(m=1\), then \(a=N\) and Corollary~\ref{prop:power-law-bond-embedding} yields the conclusion.
  If \(m\geq2\), Corollary~~\ref{prop:power-law-fractional-embedding} yields the conclusion under condition~\eqref{eq:uniform-interior-thickness}. 
  The same conclusion holds for arbitrary bounded open sets by computing a lower bound on \(\mathrm{d}\kappa_{\bx}(\boldsymbol{h})\), which 
  is needed for the application of Theorem~\ref{thm:fractional-Sobolev-embedding}, directly 
  using the product structure of Lebesgue measure in \(\boldsymbol{H}\).
  Indeed, fixing one column and integrating over the remaining \(m-1\) columns gives 
  \[ \mathrm{d}\kappa_{\bx,i}(\boldsymbol{h}) = |\O|^{m-1} \chi_{\O-\bx}(\boldsymbol{h})\,\mathrm{d}\boldsymbol{h}. \qedhere\]  
\end{proof}

%% file: section_existence.tex
%=====================================================
\section{Existence of solutions}\label{sec:existence}

\subsection{Weak lower semicontinuity under averaged convexity} 
\label{subsec:weak-lsc-averaged-convexity}

To limit the scope of this section, we restrict our attention to nonnegative energy densities.
Thus, we consider~\eqref{general} under an additional assumption \(W\geq 0\).
It will be convenient to represent the integral in~\eqref{general} in a different fashion.
Specifically, let
\[
  \boldsymbol{\mathsf{x}} = (\bx_0,\ldots,\bx_m)\in\O^{m+1}, \qquad
  \boldsymbol{\mathsf{y}} = (\boldsymbol{y}_0,\ldots,\boldsymbol{y}_m)\in(\R^d)^{m+1}.
\]
We will write
\begin{equation}\label{eq:H_G_definition}
  \boldsymbol{H}(\boldsymbol{\mathsf{x}}) = [\bx_1-\bx_0,\ldots,\bx_m-\bx_0] \in \R^{N\times m},
  \qquad
  \boldsymbol{G}(\boldsymbol{\mathsf{x}},\boldsymbol{\mathsf{y}}) = [\boldsymbol{y}_1-\boldsymbol{y}_0,\ldots,\boldsymbol{y}_m-\boldsymbol{y}_0]\boldsymbol{H}(\boldsymbol{\mathsf{x}})^\dagger \in \R^{d\times N},
\end{equation}
and
\begin{equation}\label{eq:full-particle-integrand}
  w(\boldsymbol{\mathsf{x}},\boldsymbol{\mathsf{y}}) = W(\bx_0,\boldsymbol{H}(\boldsymbol{\mathsf{x}}),\boldsymbol{y}_0,\boldsymbol{G}(\boldsymbol{\mathsf{x}},\boldsymbol{\mathsf{y}})),
\end{equation}
so that
\begin{equation}\label{eq:full-particle-representation}
  E(\bu)= \int_{\O^{m+1}} w(\boldsymbol{\mathsf{x}},\bu(\bx_0),\ldots,\bu(\bx_m))\,\mathrm{d}\boldsymbol{\mathsf{x}}.
\end{equation}
The appropriate convexity condition is imposed after averaging over the positions and values of all but one of the interacting particles. 
This is a multipoint analogue of the nonlocal convexity conditions considered in \cite{BCM14,Elbau2011}.

\begin{definition}[Particle-wise averaged convexity] 
  \label{def:particle-wise-averaged-convexity} 
  Let \(1\leq p<\infty\).
  We say that the integrand \(w\) satisfies particle-wise averaged convexity if, 
  for every \(k\in\{0,\ldots,m\}\), every collection 
  \(\bv_i\in L^p(\O;\R^d)\), \(i\neq k\), and almost every \(\bx_k\in\O\), the function 
  \begin{equation} 
    \begin{aligned} 
      \boldsymbol{y}\longmapsto \Phi_{k,\bx_k,\boldsymbol{\mathsf{v}}^{(k)}}(\boldsymbol{y}) =
      \int_{\O^m} w\bigl( \boldsymbol{\mathsf{x}}, \bv_0(\bx_0),\ldots, \bv_{k-1}(\bx_{k-1}), \boldsymbol{y}, \bv_{k+1}(\bx_{k+1}),\ldots, \bv_m(\bx_m) \bigr) 
      \prod_{i\neq k}\mathrm{d}\bx_i 
    \end{aligned} 
    \label{eq:particle-wise-averaged-convexity} 
  \end{equation} 
  is convex on \(\R^d\), where \(\boldsymbol{\mathsf{v}}^{(k)} = (\bv_i)_{i\neq k}. \)
\end{definition}

Pointwise separate convexity of \(\boldsymbol{\mathsf{y}} \longmapsto w(\boldsymbol{\mathsf{x}}, \boldsymbol{\mathsf{y}} ) \) in each column of $\boldsymbol{\mathsf{y}}$, for almost every \(\boldsymbol{\mathsf{x}}\in\O^{m+1}\), is sufficient for particle-wise averaged convexity, but it is generally stronger than necessary.  The distinction between particle-wise averaged convexity and pointwise separate convexity is clear at the formal level, although mechanically natural examples separating the two conditions are not easy to identify.  Related examples and discussion in the pairwise case \(m=1\) may be found in~\cite{BCM14}.

\begin{lemma}%[Integrated Jensen inequality] 
  \label{lem:integrated-Jensen} 
  Let \(1\leq p<\infty\), and suppose that 
  \[ w: \Omega^{m+1}\times(\mathbb R^d)^{m+1} \longrightarrow[0,+\infty] \]
  is a nonnegative normal integrand satisfying particle-wise averaged convexity 
  in the sense of Definition~\ref{def:particle-wise-averaged-convexity}. 
  For \(i=0,\ldots,m\), let 
  \(\bx_i\longmapsto\lambda_{i,\bx_i} \in\mathcal P(\R^d)\)
  be weakly measurable probability kernels\footnote{%
That is, we assume that 
\(\bx\longmapsto \int_{\R^d}\varphi(y)\,\mathrm{d}\lambda_{\bx}(\boldsymbol{y})\) 
is measurable for every \(\varphi\in C_b(\R^d)\).}
  satisfying 
  \begin{equation}\label{eq:p-moment-lambda} \int_\Omega\int_{\R^d} |\boldsymbol{y}_i|^p 
  \,\mathrm{d}\lambda_{i,\bx_i}(\boldsymbol{y}_i)\,\mathrm{d}\bx_i <+\infty. \end{equation}
  Fix \(k\in\{0,\ldots,m\}\), and let 
  \[ \overline{\lambda}_k(\bx_k) = \int_{\R^d} \boldsymbol{y}_k\,\mathrm{d}\lambda_{k,\bx_k}(\boldsymbol{y}_k). 
  \] 
  Then
  \begin{equation}\label{eq:integrated-Jensen}
  \begin{aligned} 
    \int_{\O^{m+1}} \int_{(\R^d)^m} w\bigl( \boldsymbol{\mathsf x}, \boldsymbol{y}_0,\ldots,\boldsymbol{y}_{k-1}, \overline{\lambda}_k(\bx_k), \boldsymbol{y}_{k+1},\ldots,\boldsymbol{y}_m \bigr) &\prod_{i\neq k} \mathrm{d}\lambda_{i,\bx_i}(\boldsymbol{y}_i) \,\mathrm{d}\boldsymbol{\mathsf x}\\ \leq 
    \int_{\O^{m+1}} \int_{(\R^d)^{m+1}} w\bigl( \boldsymbol{\mathsf x}, \boldsymbol{\mathsf y} \bigr) &\prod_{i=0}^{m} \mathrm{d}\lambda_{i,\bx_i}(\boldsymbol{y}_i) \,\mathrm{d}\boldsymbol{\mathsf x}.
  \end{aligned}
  \end{equation} 
\end{lemma}

\begin{proof} 
  By the randomization lemma for probability kernels~\cite[Lemma~4.22]{Kallenberg2021}, 
  for each \(i=0,\ldots,m\) there exists a measurable map 
  \( \bv_i:\O\times(0,1)\to\R^d \)
  such that 
  \( \lambda_{i,\bx} = \bv_i(\bx,\cdot)_\# \bigl( \mathcal L^1\mres(0,1) \bigr) \)
  for almost every \(\bx\in\O\).
  In particular, 
  \[ \overline{\lambda}_k(\bx_k) = \int_0^1\bv_k(\bx_k,t_k)\,\mathrm{d}t_k. \]
  The moment assumption~\eqref{eq:p-moment-lambda} and Fubini's theorem imply that 
  \(\bv_i(\cdot,t_i) \in L^p(\O;\mathbb R^d) \) for almost every \(t_i\in(0,1)\). 
  Thus, for almost every fixed \( \boldsymbol{\mathsf t}^{(k)} = (t_i)_{i\neq k}\in(0,1)^m, \) Definition~\ref{def:particle-wise-averaged-convexity} implies that the function 
  \[ \begin{aligned} 
    \boldsymbol{y}&\longmapsto \widetilde{\Phi}_{ k,\bx_k,\boldsymbol{\mathsf t}^{(k)} }(\boldsymbol{y})\\
     &= \int_{\O^m} w\Bigl( \boldsymbol{\mathsf x}, \bv_0(\bx_0,t_0),\ldots, \bv_{k-1}(\bx_{k-1},t_{k-1}), \boldsymbol{y}, \bv_{k+1}(\bx_{k+1},t_{k+1}),\ldots, \bv_m(\bx_m,t_m) \Bigr) \prod_{i\neq k}\mathrm{d}\bx_i 
  \end{aligned} \] 
  is convex for almost every \(\bx_k\in\O\). 
  For almost every \((\bx_k,\boldsymbol{\mathsf{t}}^{(k)})\), this function is also lower semicontinuous.
  Indeed, if \(\boldsymbol{y}_j\to\boldsymbol{y}\), the lower semicontinuity of 
  \(\boldsymbol{\mathsf{y}}\mapsto w(\boldsymbol{\mathsf{x}},\boldsymbol{\mathsf{y}})\) 
  and Fatou's lemma yield 
  \( \widetilde{\Phi}_{ k,\bx_k,\boldsymbol{\mathsf{t}}^{(k)} }(\boldsymbol{y}) \leq \liminf_{j\to\infty} \widetilde{\Phi}_{ k,\bx_k,\boldsymbol{\mathsf{t}}^{(k)} }(\boldsymbol{y}_j). \)

  Thus, for almost every  fixed \(\boldsymbol{\mathsf t}^{(k)}\), Jensen's 
  inequality yields 
  \[ \begin{aligned} 
    \widetilde{\Phi}_{ k,\bx_k,\boldsymbol{\mathsf t}^{(k)} } \left( \int_0^1\bv_k(\bx_k,t_k)\,\mathrm{d}t_k \right) \leq 
    \int_0^1 \widetilde{\Phi}_{ k,\bx_k,\boldsymbol{\mathsf t}^{(k)} } \bigl( \bv_k(\bx_k,t_k) \bigr) \,\mathrm{d}t_k 
   \end{aligned} \] 
   for almost every \(\bx_k\in\O\). We integrate this inequality first with respect to \(\bx_k\) and then with respect to \(\boldsymbol{\mathsf t}^{(k)}\). 
   Since \(w\geq0\), Tonelli's theorem permits all integrations to be rearranged. 
   Using the randomization identities for the kernels \(\lambda_{i,\bx_i}\) yields 
   precisely~\eqref{eq:integrated-Jensen}.
\end{proof}

This lemma is the basic property that ensures weak lower semicontinuity. 

\begin{theorem}[Weak lower semicontinuity in \(L^p(\O;\R^d)\)] \label{thm:particle-wise-averaged-lsc} 
  Let \(1\leq p<\infty\), and suppose that the full-particle integrand \(w\) defined by~\eqref{eq:full-particle-integrand} is a nonnegative normal integrand.
  Assume, in addition, that \(w\) satisfies particle-wise averaged convexity in the sense of Definition~\ref{def:particle-wise-averaged-convexity}. 
  Then \(E\) in \eqref{eq:full-particle-representation}
  is sequentially weakly lower semicontinuous on \(L^p(\O;\R^d)\).
\end{theorem}
\begin{proof} 
  Assume that
  \[ \bu_j\rightharpoonup\bu \qquad\text{in }L^p(\O;\R^d). \] 
  By switching to a subsequence, which is not relabelled, we can also assume that
  \[ \lim_{j\to\infty} E(\bu_j) = \liminf_{j\to\infty} E(\bu_j). \]
  Since \((\bu_j)\) is bounded in \(L^p(\O;\R^d)\), after passing to a further 
  subsequence it generates a Young measure 
  \[ \boldsymbol{\nu}=(\nu_{\bx})_{\bx\in\O}, \]
  see, for example, \cite{Ball1989YoungMeasures,Rindler2018}. 
  This Young measure has finite integrated \(p\)-th moment and, because \(\bu_j\rightharpoonup\bu\), its barycenter satisfies
  \begin{equation} 
    \label{eq:barycenter-identity}
    \bu(\bx) = \int_{\R^d}\boldsymbol{y}\,\mathrm{d}\nu_{\bx}(\boldsymbol{y}) 
    \qquad \text{for almost every }\bx\in\O.
  \end{equation}
  For \(p=1\), this follows from the uniform integrability implied by weak convergence in 
  \(L^1(\O)\).
  Let 
  \[
    \boldsymbol{U}_j(\bx_0,\ldots,\bx_m) = \bigl( \bu_j(\bx_0),\ldots,\bu_j(\bx_m) \bigr).
  \]
  The product property for Young measures implies that the sequence \((\boldsymbol{U}_j)\) generates the Young measure
  \(\nu_{\bx_0}\otimes\cdots\otimes\nu_{\bx_m}\)
  for almost every \(\boldsymbol{\mathsf{x}}=(\bx_0,\ldots,\bx_m)\in\O^{m+1}\), see~\cite{Pedregal1997}.
  Since \(w\) is a nonnegative normal integrand, the lower-bound theorem for Young measures~\cite[Chapter 4]{Rindler2018} and~\eqref{eq:full-particle-representation} give 
  \begin{equation} 
    \int_{\O^{m+1}} \int_{(\R^d)^{m+1}} w(\boldsymbol{\mathsf{x}},\boldsymbol{\mathsf{y}}) \prod_{i=0}^m \mathrm{d}\nu_{\bx_i}(\boldsymbol{y}_i) \,\mathrm{d}\boldsymbol{\mathsf{x}} \leq \liminf_{j\to\infty}  E(\bu_j).
     \label{eq:young-measure-lower-bound} 
  \end{equation} 
  We next compare the left-hand side of~\eqref{eq:young-measure-lower-bound} with 
  \(E(\bu)\). 

  For \(r\in\{0,\ldots,m+1\}\), let 
  \[ \lambda^{(r)}_{i,\bx} = 
      \begin{cases} \nu_{\bx}, & i<r,\\[1mm]
                    \delta_{\bu(\bx)}, & i\geq r.
      \end{cases} 
  \] 
  Their integrated \(p\)-th moments are finite, and their barycenters equal \(\bu(\bx)\),
  for almost every \(\bx\in\O\), by~\eqref{eq:barycenter-identity}.
  Let also
  \begin{equation} 
    \begin{aligned} 
        I_r &= \int_{\O^{m+1}} \int_{(\R^d)^{m+1}} w(\boldsymbol{\mathsf{x}},\boldsymbol{\mathsf{y}}) \prod_{i=0}^m \mathrm{d}\lambda^{(r)}_{i,\bx_i}(\boldsymbol{y}_i) \,\mathrm{d}\boldsymbol{\mathsf{x}}. 
    \end{aligned} 
    \label{eq:iterated-Jensen-integrals} 
  \end{equation} 
  All these integrals are well defined as elements of \([0,+\infty]\), since \(w\geq0\). 
  Applying Lemma~\ref{lem:integrated-Jensen} with \(k=r\) and with the probability kernels 
  \(\lambda_{i,\bx} = \lambda^{(r+1)}_{i,\bx}\), we obtain the inequality \(I_r\leq I_{r+1}. \) 
  Consequently, 
  \begin{equation}\label{eq:iterated-Jensen-chain} I_0\leq I_1\leq\cdots\leq I_{m+1}. \end{equation}
  For \(r=0\), every probability measure in \eqref{eq:iterated-Jensen-integrals} is a 
  Dirac mass at the corresponding barycenter. 
  Consequently, 
  \begin{equation}\label{eq:thmlsc_I0}
    I_0 = \int_{\O^{m+1}} w\bigl( \boldsymbol{\mathsf{x}}, \bu(\bx_0),\ldots,\bu(\bx_m) \bigr) \,\mathrm{d}\boldsymbol{\mathsf{x}} = E(\bu). 
  \end{equation}
   For \(r=m+1\), every particle value is integrated against the corresponding Young measure, so 
   \begin{equation}\label{eq:thmlsc_Im1} I_{m+1} = \int_{\O^{m+1}} \int_{(\R^d)^{m+1}} w(\boldsymbol{\mathsf{x}},\boldsymbol{\mathsf{y}}) \prod_{i=0}^m \mathrm{d}\nu_{\bx_i}(\boldsymbol{y}_i) \,\mathrm{d}\boldsymbol{\mathsf{x}}. \end{equation}
   Combining~\eqref{eq:thmlsc_I0}, \eqref{eq:iterated-Jensen-chain}, \eqref{eq:thmlsc_Im1},
   and~\eqref{eq:young-measure-lower-bound} proves the claimed sequential 
   lower semicontinuity.
   %Finally, weak convergence in \(NW_{\nu,\mu}^{1,p}(\Omega;\mathbb R^d)\) implies weak convergence in \(L^p(\Omega;\mathbb R^d)\), since the canonical map \[ NW_{\nu,\mu}^{1,p}(\Omega;\mathbb R^d) \longrightarrow L^p(\Omega;\mathbb R^d) \] is continuous. The last assertion follows.
  \end{proof} 

  \begin{corollary}[Lower semicontinuity on weighted nonlocal spaces] \label{cor:lsc-NW} 
  Let the assumptions of Theorem~\ref{thm:particle-wise-averaged-lsc} hold.
  Then the restriction of \(E\) to \(NW_{\nu,\mu}^{1,p}(\O;\R^d)\) is sequentially weakly lower semicontinuous.
\end{corollary} 
\begin{proof} 
  We note the continuous embedding 
  \( NW_{\nu,\mu}^{1,p}(\O;\R^d) \hookrightarrow L^p(\O;\R^d), \)
  which follows immediately from the definition of the norm~\eqref{eq:weighted-NW-norm}.
  Therefore, weak convergence \(\bu_j\rightharpoonup\bu\)
  in \( NW_{\nu,\mu}^{1,p}(\O;\R^d)\) implies
  \( \bu_j\rightharpoonup\bu\) in \(L^p(\O;\R^d).\) 
  Consequently, Theorem~\ref{thm:particle-wise-averaged-lsc} gives 
  \(E(\bu) \leq \liminf_{j\to\infty}E(\bu_j). \)
\end{proof}

Rank-one convexity of the simplicial energy density provides a practical sufficient condition for particle-wise averaged convexity.

Recall that a function  \( g: \R^{d\times N}\to \R\cup\{+\infty\} \)
is rank-one convex if 
\( t\longmapsto g( \boldsymbol{G}+t\,\boldsymbol{a}\otimes\boldsymbol{b} ) \)
is convex for every 
\(\boldsymbol{G}\in\R^{d\times N}\), 
\(\boldsymbol{a}\in\R^d\), and \(\boldsymbol{b}\in\R^N\).

\begin{proposition}[Rank-one convex energies] 
  \label{prop:rank-one-convex-averaged-convexity} 
  Suppose that the energy density has the form 
  \[ 
  W(\bx,\boldsymbol{H},\boldsymbol{y},\boldsymbol{G}) = 
  \widetilde{W}(\bx,\boldsymbol{H},\boldsymbol{G}),
  \]
  where, for almost every \((\bx,\boldsymbol{H})\), the map 
  \[ \boldsymbol{G} \longmapsto \widetilde{W}(\bx,\boldsymbol{H},\boldsymbol{G}) \] 
  is rank-one convex. 
  Then the corresponding full-particle integrand 
  \[ w( \boldsymbol{\mathsf{x}}, \boldsymbol{\mathsf{y}} ) = \widetilde{W}\left(\bx_0, \boldsymbol{H}(\boldsymbol{\mathsf{x}}), \boldsymbol{G}(\boldsymbol{\mathsf{x}}, \boldsymbol{\mathsf{y}} ) \right) 
  \] 
  is separately convex in the particle values \(\boldsymbol{y}_0,\ldots,\boldsymbol{y}_m\). In particular, \(w\) satisfies particle-wise averaged convexity in the sense of Definition~\ref{def:particle-wise-averaged-convexity}. 
\end{proposition}
\begin{proof}
  Let \(\boldsymbol{\mathsf{x}}\in\O^{m+1}\) be such that 
  \(\boldsymbol{G}\longmapsto \widetilde{W}\bigl( \bx_0, \boldsymbol{H}(\boldsymbol{\mathsf {x}}), \boldsymbol{G} \bigr) \) is rank-one convex, and fix all particle values except \(\boldsymbol{y}_k\).
  For any \(\boldsymbol{z}\in\R^d\), due to~\eqref{eq:H_G_definition} we have
  \[
    \boldsymbol{G}\bigl( \boldsymbol{\mathsf x}, \boldsymbol{y}_0,\ldots,\boldsymbol{y}_k+\boldsymbol{z},\ldots,\boldsymbol{y}_m \bigr) - 
    \boldsymbol{G}\bigl( \boldsymbol{\mathsf x}, \boldsymbol{\mathsf y} \bigr) = 
    \boldsymbol{z}\otimes\boldsymbol\xi_k( \boldsymbol{\mathsf x} ), 
  \] 
  where 
  \[ 
    \boldsymbol{\xi}_k( \boldsymbol{\mathsf x} ) = \begin{cases} 
      -[\boldsymbol{H}(\boldsymbol{\mathsf x})^{\dagger}]^{\top} 
      \displaystyle\sum_{i=1}^m\boldsymbol{e}_i, & \text{if\ } k=0,\\[3mm] 
      [\boldsymbol H(\boldsymbol{\mathsf x})^{\dagger}]^{\top} \boldsymbol{e}_k, 
      & \text{if\ } k\in\{1,\ldots,m\}. 
    \end{cases}
  \] 
  Thus, with all other particle values fixed, the map 
  \( \boldsymbol{y}_k \longmapsto \boldsymbol{G}( \boldsymbol{\mathsf x}, \boldsymbol{\mathsf{y}} ) \) is affine, and the difference between the images of any two values of \(\boldsymbol{y}_k\) has rank at most one. 
  Rank-one convexity of \(\widetilde{W}(\bx_0,\boldsymbol{H}(\boldsymbol{\mathsf{x}}),\cdot)\) therefore implies that the map \( \boldsymbol{y}_k \longmapsto w( \boldsymbol{\mathsf {x}}, \boldsymbol{\mathsf{y}} ) \) is convex with all the remaining particle values fixed.
  Finally, integration over the positions of the remaining particles preserves convexity. Hence \(w\) satisfies particle-wise averaged convexity in the sense of Definition~\ref{def:particle-wise-averaged-convexity}.
\end{proof}

Proposition~\ref{prop:rank-one-convex-averaged-convexity} gives a practical sufficient 
condition for particle-wise averaged convexity. 
It is not necessary, since the averaging in Definition~\ref{def:particle-wise-averaged-convexity} may restore convexity even when the integrand is 
not separately convex in the particle values.

\begin{remark}[Polyconvex simplicial hyperelasticity]
In the hyperelastic case \(d=N\), the conclusion of Proposition~\ref{prop:rank-one-convex-averaged-convexity} remains valid for rank-one convex densities
depending on the simplicial deformation gradient \(\boldsymbol{F} = \boldsymbol{I}
+\boldsymbol{G}.\)
In particular, polyconvexity in the sense of Ball~\cite{Ball1976} is a sufficient condition
for rank-one convexity.

Consequently, Proposition~\ref{prop:rank-one-convex-averaged-convexity} applies, in 
particular, to simplicial stored-energy densities of the form~\eqref{eq:polyconvex-simplicial-density}.
For example, in dimension \(N=3\), this includes several classical hyperelastic models such 
as the compressible Mooney--Rivlin-type energy, which has the form
\[ 
    \Psi(\boldsymbol{F}) = c_1|\boldsymbol{F}|_{F}^{2} + 
    c_2|\operatorname{cof}\boldsymbol{F}|_{F}^{2} + h(\det\boldsymbol{F}),
\] 
where \(c_1,c_2\geq0\), and \(h:\R\to[0,+\infty] \) is convex and lower semicontinuous, and satisfies \( h(t)=+\infty\) for \(t\leq0.\)
The choice \(c_2=0\) includes compressible neo-Hookean-type energies. 
\end{remark}

%===================================================== 
\subsection{Existence under averaged convexity} 
\label{subsec:existence-averaged-convexity} 

We now combine the weak lower semicontinuity result of Subsection~\ref{subsec:weak-lsc-averaged-convexity} with the nonlocal Poincar{\'e} inequality established in Subsection~\ref{subsec:nonlocal-Poincare}.
Let \(\O_D\subset\O\) be measurable, with \(|\O_D|>0, \) and let 
\[ \bu_D \in NW_{\nu,\mu}^{1,p}(\O;\R^d) \] be prescribed.
We define the corresponding admissible class by
\begin{equation} 
  \begin{aligned} 
    \mathcal{A}_{\bu_D} &= \left\{ \bu\in NW_{\nu,\mu}^{1,p}(\O;\R^d) \colon 
    \bu=\bu_D \text{ almost everywhere on }\O_D \right\} \\ &= \bu_D+ NW_{\nu,\mu,D}^{1,p}(\O;\R^d). 
  \end{aligned} 
  \label{eq:nonhomogeneous-Dirichlet-class}
\end{equation} 
We interpret the condition in~\eqref{eq:nonhomogeneous-Dirichlet-class} as a nonlocal inhomogeneous Dirichlet condition. 
The admissible class is a closed affine subspace of \(NW_{\nu,\mu}^{1,p}(\O;\R^d)\).
Indeed, if \(\bu_j\to\bu\) in \(NW_{\nu,\mu}^{1,p}\), then \(\bu_j\to\bu\) in \(L^p(\O;\R^d)\), and consequently also in \(L^p(\O_D;\R^d)\).
It follows that \(\bu=\bu_D\) almost everywhere on \(\Omega_D\). 
Since a closed convex subset of a Banach space is weakly closed, \(\mathcal{A}_{\bu_D}\) is also weakly closed.
We next record the coercivity consequence of the nonlocal Poincar{\'e} inequality. 
\begin{proposition}[Coercivity on the Dirichlet class] 
  \label{prop:coercivity-Dirichlet-class} 
  Let \(1\leq p<\infty\), and suppose that the assumptions of Proposition~\ref{prop:NW-Poincare} hold for the prescribed set \(\O_D\). 
  Suppose further that there exist \(c_0>0\) and \(C_0\geq0\) such that 
  \[ 
  E(\bu) \geq c_0 [\bu]_{NW_{\nu,\mu}^{1,p}(\O;\R^d)}^p - C_0, 
  \qquad \forall \bu\in\mathcal{A}_{\bu_D},
  \] 
  where \(\mathcal{A}_{\bu_D}\) is given by~\eqref{eq:nonhomogeneous-Dirichlet-class}. 
  Then there exist \(c>0\),\(C\geq 0\), depending also on \(\bu_D\), such that 
  \[ 
  E(\bu) \geq c \|\bu\|_{NW_{\nu,\mu}^{1,p}(\Omega;\mathbb R^d)}^p - C, 
  \qquad \forall \bu\in\mathcal{A}_{\bu_D}.\] 
\end{proposition} 
\begin{proof} 
  Take an arbitrary \(\bu\in\mathcal{A}_{\bu_D}\), and set 
  \(\bv=\bu-\bu_D\) so that \(\bv \in NW_{\nu,\mu,D}^{1,p}(\O;\R^d).\) 
  By Proposition~\ref{prop:NW-Poincare}, there exists \(C_D>0\) such that \( \|\bv\|_{L^p(\O;\R^d)}^p \leq C_D [\bv]_{NW_{\nu,\mu}^{1,p}(\O;\R^d)}^p. \) 
  Moreover, \( [\bv]_{NW_{\nu,\mu}^{1,p}(\O;\R^d)} \leq [\bu]_{NW_{\nu,\mu}^{1,p}(\O;\R^d)} + [\bu_D]_{NW_{\nu,\mu}^{1,p}(\O;\R^d)}. \) 
  Consequently, 
  \[\begin{aligned} 
    \|\bu\|_{L^p(\O;\R^d)}^p 
    &\leq 2^{p-1} \left( \|\bv\|_{L^p(\O;\R^d)}^p+\|\bu_D\|_{L^p(\O;\R^d)}^p \right) \\ 
    &\leq 2^{p-1}C_D [\bv]_{NW_{\nu,\mu}^{1,p}(\O;\R^d)}^p + 2^{p-1}\|\bu_D\|_{L^p(\O;\R^d)}^p \\ 
    &\leq 2^{2p-2}C_D [\bu]_{NW_{\nu,\mu}^{1,p}(\O;\R^d)}^p + C_1, 
  \end{aligned}\] 
  where \(C_1\geq 0\) depends on \(p\), \(C_D\), and \(\bu_D\). 
  Thus, there exist \(C_2>0\) and \(C_3\geq0\) such that \(\|\bu\|_{NW_{\nu,\mu}^{1,p}(\O;\R^d)}^p \leq C_2 [\bu]_{NW_{\nu,\mu}^{1,p}(\O;\R^d)}^p + C_3. \) 
  Combining this inequality with the assumed lower bound for \(E\) yields the conclusion.
\end{proof}

Combining the preceding coercivity result with sequential weak lower semicontinuity yields the following abstract existence theorem, which amounts to the application of the direct method.

\begin{theorem}
  \label{thm:existence-direct-method} 
  Let \(1<p<\infty\), and assume that: 
  \begin{enumerate} 
    \item the functional \(E: NW_{\nu,\mu}^{1,p}(\O;\R^d) \longrightarrow \R\cup\{+\infty\} \) is
    sequentially weakly lower semicontinuous; 
    \item the assumptions of Proposition~\ref{prop:NW-Poincare} (nonlocal Poincar{\'e} inequality) hold for the prescribed set \(\O_D\); 
    \item there exist \(c_0>0\) and \(C_0\geq0\) such that 
    \[ E(\bu) \geq c_0 [\bu]_{NW_{\nu,\mu}^{1,p}(\O;\R^d)}^p - C_0,
    \qquad \forall \bu\in\mathcal A_{\bu_D},\]
    where \(\mathcal{A}_{\bu_D}\) is given by~\eqref{eq:nonhomogeneous-Dirichlet-class};
    \item there exists \(\bu_0\in\mathcal{A}_{\bu_D}\) such that 
    \(E(\bu_0)<+\infty.\)
  \end{enumerate} 
  Then the infimum \(\alpha = \inf_{\bu\in\mathcal{A}_{\bu_D}} E(\bu)\) is attained.
\end{theorem} 
\begin{proof}
  By the coercivity hypothesis, \(\alpha>-\infty, \) while the existence of 
  \(\bu_0\in\mathcal A_{\bu_D}\) with finite energy yields \(\alpha\leq E(\bu_0)<+\infty. \)
  Let \((\bu_j)\subset\mathcal{A}_{\bu_D}\) be a minimizing sequence.
  After discarding finitely  many terms, the sequence  \((E(\bu_j))\) is bounded. 
  By Proposition~\ref{prop:coercivity-Dirichlet-class}, \((\bu_j)\) is bounded in \(NW_{\nu,\mu}^{1,p}(\O;\R^d)\). 
  Since \(1<p<\infty\), \(NW_{\nu,\mu}^{1,p}(\O;\R^d)\) is reflexive. 
  Consequently, after passing to a further subsequence, there exists 
  \( \bar{\bu}\in NW_{\nu,\mu}^{1,p}(\O;\R^d) \) such that \( \bu_j\rightharpoonup\bar{\bu}\) in \(NW_{\nu,\mu}^{1,p}(\O;\R^d).\)
  The admissible class \(\mathcal{A}_{\bu_D}\) is weakly closed, and hence 
  \( \bar{\bu}\in\mathcal{A}_{\bu_D}\) and \(\alpha\leq E(\bar{\bu})\).
  Finally, sequential weak lower semicontinuity gives 
  \( E(\bar{\bu}) \leq \liminf_{j\to\infty} E(\bu_j) = \alpha \),
  which shows that the infimum is attained at \(\bar{\bu}\).
\end{proof}

\begin{corollary}[Existence under averaged convexity] 
  \label{cor:existence-averaged-convexity} 
  Let \(1<p<\infty\), and suppose that the assumptions of 
  Theorem~\ref{thm:particle-wise-averaged-lsc} (averaged convexity),  
  and Proposition~\ref{prop:coercivity-Dirichlet-class} (coercivity) hold. 
  Assume, in addition, that there exists \(\bu_0\in\mathcal{A}_{\bu_D}\) such that 
  \(E(\bu_0)<+\infty.\)
  Then \(E\) attains its infimum on \(\mathcal{A}_{\bu_D}\). 
\end{corollary} 
\begin{proof} 
  Due to Corollary~\ref{cor:lsc-NW}, the functional \(E\) is 
  sequentially weakly lower semicontinuous on \(NW_{\nu,\mu}^{1,p}(\O;\R^d)\).
  The conclusion therefore follows from Theorem~\ref{thm:existence-direct-method}. 
\end{proof}

\begin{remark} 
  The seminorm-coercivity assumption in Proposition~\ref{prop:coercivity-Dirichlet-class}
  follows directly from a pointwise \(p\)-growth condition when the interaction measure in 
  the energy agrees with those defining the weighted nonlocal space. 
  For example, if 
  \[ \begin{aligned} 
    E_{\nu,\mu}(\bu) &= \int_\O \int_{\mathbb{S}^{Nm-1}} \int_{r_\O(\bx,\boldsymbol{S})}
     W\bigl( \bx,r\boldsymbol{S},\bu(\bx), \bD\bu(\bx;r\boldsymbol{S}) \bigr) 
     \,\mathrm{d}\nu(r)\,\mathrm{d}\mu(\boldsymbol{S})\,\mathrm{d}\bx 
  \end{aligned} \] 
  and 
  \[ 
    W( \bx,r\boldsymbol{S},\boldsymbol{y},\boldsymbol{G} ) \geq 
    c_0|\boldsymbol{G}|_{F}^p - a(\bx,\boldsymbol{S},r), 
  \] 
  where \(a\) is integrable with respect to the underlying interaction measure, then 
  \[ 
    E_{\nu,\mu}(\bu) \geq c_0[\bu]_{NW_{\nu,\mu}^{1,p}(\O;\R^d)}^p - \|a\|_{L^1}. 
  \] 
  For the energy~\eqref{general}, integrated with respect to Lebesgue measure 
  in \(\boldsymbol{H}\), the corresponding growth condition directly controls the 
  unweighted \(NW^{1,p}(\O;\R^d)\)-seminorm. 
  A comparison of interaction measures is required to obtain coercivity in a different
  weighted space.
\end{remark}

%===================================================== 
\subsection{Existence by fractional compactness} 
\label{subsec:existence-fractional-compactness} 
The averaged-convexity condition in Subsection~\ref{subsec:weak-lsc-averaged-convexity} 
provides weak lower semicontinuity without requiring compactness.
Convexity can instead be dispensed with when the nonlocal space is compactly embedded 
into \(L^p\).
Indeed, strong \(L^p\)-convergence gives, after passage to a subsequence, pointwise
convergence of all the particle values, and lower semicontinuity then follows directly 
from Fatou's lemma. 

\begin{proposition}[Lower semicontinuity by compactness] 
  \label{prop:lsc-by-compactness} 
  Let \(1\leq p<\infty\), and suppose that 
  \(NW_{\nu,\mu}^{1,p}(\O;\R^d) \hookrightarrow L^p(\O;\R^d)\) compactly. 
  Assume that the full-particle integrand 
  \(w: \O^{m+1}\times(\R^d)^{m+1} \longrightarrow [0,+\infty]\)
  defined by~\eqref{eq:full-particle-integrand} is a nonnegative normal integrand. 
  Then \(E\) is sequentially weakly lower semicontinuous on 
  \(NW_{\nu,\mu}^{1,p}(\O;\ R^d)\), and no convexity assumption on \(w\) is required. 
\end{proposition} 
\begin{proof} 
  Let \(\bu_j\rightharpoonup\bu\) in \(NW_{\nu,\mu}^{1,p}(\O;\R^d). \)
  Passing to a subsequence, which is not relabelled, we may assume that 
  \(\lim_{j\to\infty}E(\bu_j) = \liminf_{j\to\infty}E(\bu_j). \) 
  The sequence \((\bu_j)\) is bounded in \(NW_{\nu,\mu}^{1,p}(\O;\R^d)\). 
  By the assumed compact embedding, after passing to a further subsequence, there exists 
  \(\bv\in L^p(\O;\R^d)\) such that \(\bu_j\longrightarrow\bv\), strongly in \(L^p(\O;\R^d)\).
  On the other hand, the continuous embedding \(NW_{\nu,\mu}^{1,p}(\O;\R^d) \hookrightarrow L^p(\O;\R^d) \) implies that \(\bu_j\rightharpoonup\bu\), in \(L^p(\O;\R^d)\). 
  Therefore \(\bv=\bu\), and hence \(\bu_j\longrightarrow\bu\), strongly in \(L^p(\O;\R^d)\). 
  Passing to another subsequence if necessary, we obtain 
  \(\bu_j(\bx)\longrightarrow\bu(\bx),\) for almost every \(\bx\in\O\).
  Consequently, 
  \[ 
  \bigl( \bu_j(\bx_0),\ldots,\bu_j(\bx_m) \bigr) \longrightarrow 
  \bigl( \bu(\bx_0),\ldots,\bu(\bx_m) \bigr),
  \qquad  \text{for almost every\ } 
  \boldsymbol{\mathsf{x}} = (\bx_0,\ldots,\bx_m) \in\O^{m+1}. 
  \] 
  Since \(w\) is a normal integrand, the map \(\boldsymbol{\mathsf{y}} \longmapsto w( \boldsymbol{\mathsf{x}}, \boldsymbol{\mathsf{y}} ) \) is lower semicontinuous, for almost every \(\boldsymbol{\mathsf{x}}\in\O^{m+1}\). 
  Therefore, 
  \[ 
    w\left( \boldsymbol{\mathsf x}, \bu(\bx_0),\ldots,\bu(\bx_m) \right)
    \leq \liminf_{j\to\infty} w\left( \boldsymbol{\mathsf{x}}, \bu_j(\bx_0),\ldots,\bu_j(\bx_m) \right),
  \]
  for almost every \(\boldsymbol{\mathsf{x}}\in\O^{m+1}\). 
  Fatou's lemma, and~\eqref{eq:full-particle-representation} now yields 
  \[ \begin{aligned} 
    E(\bu) &= 
    \int_{\O^{m+1}} w\left( \boldsymbol{\mathsf x}, \bu(\bx_0),\ldots,\bu(\bx_m) \right) \,\mathrm{d}\boldsymbol{\mathsf{x}} 
    \leq \liminf_{j\to\infty} 
    \int_{\O^{m+1}} w\left( \boldsymbol{\mathsf x}, \bu_j(\bx_0),\ldots,\bu_j(\bx_m) \right) \,\mathrm{d}\boldsymbol{\mathsf{x}} 
    \\&= \liminf_{j\to\infty}E(\bu_j). 
  \end{aligned} \] 
  This proves the claimed sequential weak lower semicontinuity. 
\end{proof}

\begin{corollary}[Existence by compactness] 
  \label{cor:existence-by-compactness} 
  Let \(1\leq p<\infty\), and suppose that: 
  \begin{enumerate} 
    \item the embedding \(NW_{\nu,\mu}^{1,p}(\O;\R^d) \hookrightarrow L^p(\O;\R^d) \) is compact; 
    \item the full-particle integrand \(w\) is a nonnegative normal integrand; 
    \item the assumptions of Proposition~\ref{prop:coercivity-Dirichlet-class} hold; 
    \item there exists \(\bu_0\in\mathcal{A}_{\bu_D}\) such that \(E(\bu_0)<+\infty. \) 
  \end{enumerate} 
  Then \(E\) attains its infimum on \(\mathcal{A}_{\bu_D}\). 
  No convexity assumption on \(w\) is required. 
\end{corollary} 
\begin{proof} 
  Suppose first that \(1<p<\infty\).
  By Proposition~\ref{prop:lsc-by-compactness}, the functional \(E\) is sequentially 
  weakly lower semicontinuous on \(NW_{\nu,\mu}^{1,p}(\O;\R^d)\). 
  The conclusion follows from Theorem~\ref{thm:existence-direct-method}. 
  
  It remains to consider \(p=1\). 
  Let \(\alpha = \inf_{\bu\in\mathcal{A}_{\bu_D}}E(\bu)\).
  The existence of a finite-energy competitor \(\bu_0\) and the coercivity estimate imply that \( -\infty<\alpha<+\infty. \)  
  Let 
  \((\bu_j)\subset\mathcal{A}_{\bu_D}\) be a minimizing sequence.
  By Proposition~\ref{prop:coercivity-Dirichlet-class}, \((\bu_j)\) is bounded in 
  \(NW_{\nu,\mu}^{1,1}(\O;\R^d)\). 
  The compact embedding therefore gives, after passing to a subsequence, \( \bu_j\longrightarrow\bu\),  strongly in \(L^1(\O;\R^d)\) for some \(\bu\in L^1(\O; \R^d)\). 
  Passing to a further subsequence, we may assume that \(\bu_j(\bx)\longrightarrow\bu(\bx)\), for almost every \(\bx\in\O.\)
  For every fixed admissible \((\boldsymbol S,r)\), translations preserve Lebesgue-null sets. 
  Fubini's theorem therefore gives
  \(\Delta_{r\boldsymbol{S}}\bu_j(\bx) \longrightarrow \Delta_{r\boldsymbol{S}}\bu(\bx), \) and consequently \(\bD\bu_j(\bx;r\boldsymbol{S}) \longrightarrow \bD\bu(\bx;r\boldsymbol{S}), \)  for almost every admissible \((\bx,\boldsymbol S,r)\).
  Fatou's lemma gives 
  \[[\bu]_{NW_{\nu,\mu}^{1,1}(\O;\R^d)} \leq \liminf_{j\to\infty} [\bu_j]_{NW_{\nu,\mu}^{1,1}(\O;\R^d)} <+\infty. \]
  Thus \(\bu\in NW_{\nu,\mu}^{1,1}(\O;\R^d). \)
  Since \(\bu_j=\bu_D\), almost everywhere on \(\O_D\), and \(\bu_j\to\bu\) in \(L^1(\O_D;\R^d)\), we have \(\bu=\bu_D\), almost everywhere on \(\O_D. \) 
  Therefore \(\bu\in\mathcal{A}_{\bu_D}\). 
  Finally, 
  \[ \bigl( \bu_j(\bx_0),\ldots,\bu_j(\bx_m) \bigr) \longrightarrow \bigl( \bu(\bx_0),\ldots,\bu(\bx_m) \bigr),
  \qquad \text{for almost every\ } (\bx_0,\ldots,\bx_m)\in\O^{m+1}.\]
  The lower semicontinuity of \(w(\boldsymbol{\mathsf{x}},\cdot)\) and Fatou's lemma yield 
  \[ E(\bu) \leq \liminf_{j\to\infty}E(\bu_j) = \alpha. \] 
  Hence \(E(\bu)=\alpha\), and the infimum is attained.
\end{proof}

\begin{corollary}[Existence by fractional compactness] 
  \label{cor:existence-fractional-compactness} 
  Let \(1\leq p<\infty\), and suppose that the assumptions of Theorem~\ref{thm:fractional-Sobolev-embedding} hold for some \(s\in(0,1)\). 
  Assume, in addition, that \(W^{s,p}(\O;\R^d) \hookrightarrow L^p(\O;\R^d) \) compactly,
  as is the case, for example, when \(\O\) is a bounded Lipschitz domain. 
  Suppose further that: 
  \begin{enumerate} 
    \item the full-particle integrand \(w\) is a nonnegative normal integrand; 
    \item the assumptions of Proposition~\ref{prop:coercivity-Dirichlet-class} hold; 
    \item there exists \(\bu_0\in\mathcal{A}_{\bu_D}\) such that \(E(\bu_0)<+\infty. \) 
  \end{enumerate} 
  Then \(E\) attains its infimum on \(\mathcal{A}_{\bu_D}\), without any convexity assumption on \(w\). 
\end{corollary} 
\begin{proof} 
  Theorem~\ref{thm:fractional-Sobolev-embedding} gives the continuous embedding 
  \(NW_{\nu,\mu}^{1,p}(\O;\R^d) \hookrightarrow W^{s,p}(\O;\R^d). \) 
  Composing it with the compact embedding \(W^{s,p}(\O;\R^d) \hookrightarrow L^p(\O;\R^d) \) shows that \(NW_{\nu,\mu}^{1,p}(\O;\R^d) \hookrightarrow L^p(\O;\R^d) \) compactly. 
  The conclusion follows from Corollary~\ref{cor:existence-by-compactness}. 
\end{proof}

The lower-semicontinuity mechanisms in Subsections~\ref{subsec:weak-lsc-averaged-convexity} 
and \ref{subsec:existence-fractional-compactness} are complementary. 
Particle-wise averaged convexity furnishes lower semicontinuity under weak \(L^p\)-convergence
without requiring compactness of the nonlocal space. 
By contrast, fractional compactness yields strong \(L^p\)-convergence and removes the
convexity requirement entirely, at the cost of stronger assumptions on the edge marginals.